\documentclass[11pt,preqno,a4paper]{amsart}

\usepackage{enumitem}
\usepackage{amssymb,amsmath,mathtools}
\usepackage{hyperref}
\usepackage{amsmath}
\usepackage[nobysame]{amsrefs}
\usepackage{xcolor}

\usepackage{microtype}
\usepackage{bbm}
\usepackage{a4wide}
\newcommand*{\mailto}[1]{\href{mailto:#1}{\nolinkurl{#1}}}

 \allowdisplaybreaks

\newtheorem{theorem}{Theorem}[section]
\newtheorem{assumption}{Assumption}[section]
\newtheorem{lemma}[theorem]{Lemma}
\newtheorem{definition}[theorem]{Definition}
\newtheorem{proposition}[theorem]{Proposition}
\newtheorem{corollary}[theorem]{Corollary}
\newtheorem{remark}[theorem]{Remark}

\newcommand{\R}{{\mathbb R}}
\newcommand{\F}{{\mathcal F}}
\newcommand{\N}{{\mathbb N}}

\newcommand{\be}{\begin{equation}}
\newcommand{\ee}{\end{equation}}

\newcommand{\eq}{\begin{equation}}
\newcommand{\eeq}{\end{equation}}

\numberwithin{equation}{section}

\makeatletter
\newcommand{\dlmf}[1]{%
	\cite[%
	\def\nextitem{\def\nextitem{, }}%
	\@for \el:=#1\do{\nextitem\href{http://dlmf.nist.gov/\el}{(\el)}}%
	]{dlmf}%
}
\makeatother

\begin{document}

\allowdisplaybreaks

\title[Decomposition of global solutions of bi-Laplacian Nonautonomous  Schr\"odinger equations ]{Decomposition of global solutions  of bi-laplacian Nonautonomous Schr\"odinger equations}

 \author[Avy Soffer]{Avy Soffer}
 \address[Avy Soffer]{\newline
        Department of Mathematics, \newline
         Rutgers University, New Brunswick, NJ 08903 USA.}
  \email[]{soffer@math.rutgers.edu}

 \author[Jiayan Wu]{Jiayan Wu}
 \address[Jiayan Wu]{\newline
        School of mathematics, South China University of Technology, \\
 Guangzhou 510641, China.}
  \email[]{wujiayan@scut.edu.cn}

\author[Xiaoxu Wu]{Xiaoxu Wu}
\address[Xiaoxu Wu]{\newline
       Mathematical Sciences Institute, \newline Australia National University, Canberra 2601, Australia.}
 \email[]{xiaoxu.wu@anu.edu.au}

\author[Ting Zhang]{Ting Zhang}
\address[Ting Zhang]{\newline
        School of mathematics, Zhejiang university, \\
 Hangzhou 310058, China.}
 \email[]{zhangting79@zju.edu.cn}



\date{\today}


\begin{abstract} 
We study the bi-Laplacian Schr\"odinger equation with a general interaction term, which may be linear or nonlinear and is allowed to be time-dependent.  
We show that global solutions to such equations decompose asymptotically into a free wave and a weakly localized component in all space dimensions.  
Moreover, in dimensions $n \geq 9$, we prove that the weakly localized component is in fact spatially localized.  
The proof is based on a suitably adapted construction of the Free Channel Wave Operator, building on the method recently developed in~\cite{SW20221}.

\end{abstract}

\maketitle

\section{Introduction}
\subsection{Background}
{The analysis of scattering—and in particular the proof of Asymptotic Completeness (AC)—goes back to the 1920s. The rigorous literature initially focused on time-independent perturbations of the Laplacian that are localized or small \cites{MR2038194,MR2094850}. Functional-analytic techniques were used to study the spectral properties of such operators. This line of work was later extended to $N$-body Hamiltonians, beginning with the classical papers of Faddeev on the three-body problem in 1963 \cites{MR0163695,MR0160573,MR0149843}.
These methods, however, did not apply to time-dependent potentials or to nonlinear equations. In 1978, Enss \cite{Enss1978} introduced a new, time-dependent approach that opened new avenues for scattering theory, avoiding the direct analysis of the Green’s function of the Hamiltonian. Building on this idea led to the spectral methods developed by Mourre (1979–1983) \cites{MR0723552,MR0603501,MR0539739}. These developments also yielded new, comprehensive proofs of AC for three-body short- and long-range scattering. Between 1986 and 1989, Sigal and Soffer \cites{MR0818063,MR0921295,MR0915575,MR0898052,MR1034551} introduced another approach: using general estimates arising from Mourre’s method, combined with phase-space analysis, they derived a-priori propagation estimates for the full dynamics and proved AC for $N$-body Hamiltonians, along with other propagation results of interest. Nevertheless, major progress on genuinely time-dependent perturbations remained elusive.

For nonlinear evolutions whose solutions remain spatially localized for all times, we expect the solution to evolve as a superposition of a linear solution and coherent structures—such as solitons, breathers, kinks, black holes, and vortices. In many cases, we believe these coherent structures coincide with solitons of the corresponding autonomous nonlinear equations. This is the 'soliton resolution conjecture'. It asserts that if $\psi(t)$ is a global solution to the Nonlinear Schr\"odinger equation (NLS), then as $t\to\pm\infty$ the solution decomposes into a free wave $\psi_{\mathrm{free}}(t)$ and a (finite sum of) solitons $\psi_{\mathrm{sol}}(t)$:

$$
\lim_{t\to\pm\infty}\bigl\|\psi(t)-\psi_{\mathrm{free}}(t)-\psi_{\mathrm{sol}}(t)\bigr\|_{L_x^2(\mathbb{R}^n)}=0 .
$$


Between 2004 and 2014, Tao \cites{Tao1,Tao2,Tao3,T2014} developed foundational estimates for NLS with large initial data and in multichannel regimes (where solitons occur). These results enabled, for the first time, the resolution of broad classes of scattering problems without relying on the a priori decay estimates typically derived from spectral theory. Although obtained in special settings, they opened the door to substantial progress. Kenig and Merle initiated the study of the energy-critical wave equation in \cite{kenig2008global}. Building on this foundation, a sequence of works—including major collaborations with Collot and Duyckaerts \cites{duyckaerts2014scattering, collot2024soliton, duyckaerts2019soliton} for dimensions 6, 3, and all odd spatial dimensions, as well as the references cited in \cite{collot2024soliton}—pushed the theory forward. Complementing these advances, Jendrej and Lawrie \cite{JL2023} provided key insights in dimensions 4 and higher, playing a central role in resolving the asymptotic dynamics of the wave equation with an energy-critical interaction term. Taken together, these contributions cover all dimensions $3$ and above for spherically symmetric initial data in $\dot{H}^1$.

A more general phase-space analysis of scattering, which is suitable for nonlinear equations was introduced in \cite{Liu-Soffer2022}.
It provided the construction of the Free channel and weakly localized states for NLS type equations in the Radially symmetric case.
A more general phase-space method followed, in which for the first time, a new formalism is used for the construction of the free channel  \cites{SW20221,SW20224}.
In this work we follow an adapted approach to scattering theory, introduced by Soffer and Wu \cites{SW20221,SW20224}. This method yields a priori propagation estimates for broad classes of dispersive and hyperbolic equations, both time-dependent linear and nonlinear. It allows the construction of asymptotic states in great generality; in particular, it proves the existence of free-channel scattering for very general equations with localized interactions, without symmetry assumptions, and gives comprehensive results for the standard NLS and Nonlinear Klein-Gordon equations. Here we employ this framework to study the NLS with bi-Laplacian dynamics. The bi-Laplacian (bi-harmonic) Schrödinger equation was introduced by Karpman \cite{Karpman1} and Karpman–Shagalov \cite{MR1396248} to incorporate small fourth-order dispersive effects in the propagation of intense laser beams in bulk media with Kerr nonlinearity. A weak soliton-resolution result for the bi-Laplacian dynamics was obtained by Roy \cite{MR3651204} using Tao’s approach \cite{T2014}.

The approach we use here is different from Tao's. In particular, we derive localization of the weakly localized part, with explicit decay rates. This was done earlier for NLS in \cite{soffer2023soliton}. See also the review \cite{soffer2024new}.
The proof of localization is different from previous works.

\subsection{Models }Let $H_0:=(-\Delta)^2$. We consider the general class of NLSE with bi-Laplacian dynamics:
\begin{equation}\label{fourth-order}
i\partial_t\psi(t)=H_0\psi(t)+\mathcal{N}(x,t,|\psi(t)|)\psi(t),\quad (x,t)\in \R^{n+1}, n\geq 1,
\end{equation}
with initial data $\psi(t)\vert_{t=0}=\psi_0\in L^2_x(\mathbb{R}^n) $. Here we use the convention $\psi(t)\equiv \psi(x,t)$ and $\mathcal{N}$ may be either linear or nonlinear, and may be time-dependent: 
\eq
\mathcal{N}(x,t,|\psi(t)|)=V(x,t),\, N(|\psi(t)|), \text{ or } V(x,t)+N(|\psi(t)|). 
\eeq
Note that $\mathcal{N}$ need not be real-valued. Throughout this paper, we assume that \eqref{fourth-order} admits a global solution in $L^{2}(\mathbb{R}^{n})$ for all $t\in\mathbb{R}$.
\begin{assumption}\label{asp:sol}
We assume that the initial datum $\psi(x,0)=\psi_0(x)$ leads to a global solution  in $H^2\equiv H^2(\mathbb{R}^n)$ with the $L^2$ norm uniformly bounded in $t$. That is, $\psi(t)$ exists in $H^2$ for all $t\in \mathbb{R}$ and 
\begin{equation}\label{eq: uniformL2}
m:=\sup\limits_{t\in \mathbb R}\|\psi(t)\|_{L^2_x(\mathbb R^n)}<\infty.
\end{equation}
\end{assumption}
\begin{remark}Part of Assumption \ref{asp:sol} can be rephrased to state that the maximum interval of existence of $\psi(t)$ in $H^2$ is $\mathbb{R}$. 
\end{remark}
Using Assumption \ref{asp:sol}, the term $\mathcal{N}(x,t,|\psi(t)|)$ is well defined for all $t\in\mathbb{R}$, and we consider nonlinearity as a time-dependent perturbation. For notational convenience, when no confusion can arise we sometimes write $\mathcal{N}(x,t,|\psi(t)|)=V(x,t)$. Let $\langle x\rangle=(1+|x|^2)^{1/2}$, and define the weighted space,
\begin{align}
L^{p}_{x,\sigma}(\mathbb{R}^n):=\{f(x): \langle x\rangle^{\sigma}f(x)\in L^p_x(\mathbb{R}^n)\}
\end{align}
for $1\le p\le\infty$. We assume that the interaction $\mathcal{N}(x,t,|\psi(t)|)$ (equivalently, $V(x,t)$) satisfies the following condition:
\begin{assumption}\label{condition}
When $n\ge 1$, the interaction $\mathcal N(x,t,|\psi(t)|)$ and the solution $\psi(t)$ obey
\begin{equation}
\langle x\rangle^\sigma \mathcal N(x,t,|\psi(t)|),\psi(t)\in L^\infty_tL^2_x(\mathbb{R}^{n+1})
\end{equation}
for some $\sigma>1$.
\end{assumption}

\begin{assumption}\label{condition1} When $n\geq 5$, the interaction $\mathcal N(x,t,|\psi(t)|)$ satisfies 
\begin{equation}
   \mathcal N(x,t,|\psi(t)|)\in L^\infty_t L^2_x(\mathbb R^{n+1}).
\end{equation}
\end{assumption}


 

{Typical examples of Assumption \ref{condition} are $\mathcal{N}(x,t,|\psi(t)|)=\pm\lambda|\psi(t)|^p$, $p\in(1,2]$, $\lambda>0$, in three spatial dimensions under a radial symmetry assumption on the initial data.}\par
When $n=3$, if $\psi(0)$ is radial in $x$ and if
\begin{equation}
E:=\sup_{t\in \R}\|\psi(t)\|_{H^1(\R^3)}<\infty,
\end{equation}
then
\begin{align}
\sup_{t\in \R}|\psi(t)|_{L^6_x(\R^3)}\lesssim_{E} 1,
\end{align}
and
\begin{align}
|\psi(t)|\lesssim_{E}\frac{1}{|x|}, \ \text{ when } |x|\geq 1.
\end{align}
Therefore, by interpolation and the unitarity of $U(t,0)$, for $p\in(1,2]$, taking $\sigma=p$ and $q>0$ such that
\begin{equation}
\frac{1}{q}+\frac{p+1}{6}=\frac{1}{2},
\end{equation}
we obtain
\begin{align}
\|\chi(|x|<1)\langle x\rangle^{\sigma}|\psi|^p\psi\|_{L^2_x(\R^3)}\lesssim  \|\chi(|x|<1)\|_{L^q_x(\mathbb R^3)}\sup_{t\in \R}\|\psi(t)\|^{p+1}_{L^6_x(\R^3)}\lesssim_E 1
\end{align}
and
\begin{align}
\|\chi(|x|\geq1)\langle x\rangle^{\sigma}|\psi(t)|^p\psi\|_{L^2\_x(\R^3)}\lesssim_E  1.
\end{align}
Thus, $\mathcal{N}(x,t,|\psi(t)|)=\pm \lambda |\psi(t)|^p\in L^\infty_tL^2_{x,\sigma}(\R^{3+1})$.

This paper focuses on the scattering theory for \eqref{fourth-order}. Its central goal is to establish AC in the following sense:
\begin{equation}\label{AC}
\big\| \psi(t)-\textstyle\sum_{a}\psi_{a,\pm}(t)\big\|_{L^2_x(\mathbb{R}^n)}\longrightarrow 0
\quad\text{as } t\to \pm\infty,
\end{equation}
where $\{\psi_{a,\pm}\}$ are model profiles that solve simpler reference dynamics. A basic example of a
reference dynamics is the free evolution, in which case $\psi_{a,\pm}$ are \emph{free waves}.

We call an initial datum $\psi_0$ a \emph{(free) scattering state} for \eqref{fourth-order} if there exist
$\psi_\pm\in L_x^2(\mathbb{R}^n)$ such that
\[
\big\|\psi(t)-e^{-itH_0}\psi_\pm\big\|_{L_x^2(\mathbb{R}^n)}\to 0 \quad \text{as } t\to\pm\infty,
\]
i.e., the asymptotic decomposition \eqref{AC} consists purely of the free channel.

Another example of $\psi_{a,\pm}$ (or another channel) is a \emph{soliton}. Recall that the soliton-resolution conjecture is an analogue of AC for autonomous nonlinear Schr\"odinger equations; a (standing-wave) soliton $\psi_{\mathrm{sol}}(t)$
is a stable solution satisfying
\[
\big(H_0+\mathcal{N}(|\psi_{\mathrm{sol}}|)\big)\psi_{\mathrm{sol}}=\lambda\,\psi_{\mathrm{sol}}
\quad\text{for some }\lambda\in\mathbb{R}.
\]
Scattering theory studies the long-time behavior of such solutions. Note that the initial datum need not be a (free) scattering state; to try to address the soliton-resolution conjecture one must identify the free component and provide the properties of the non-free (e.g. a bound state/soliton) component as many as possible.


Proving soliton resolution or AC in the presence of solitons is challenging. In particular, the set of scattering states for nonlinear or time-dependent systems is subtle. It is known that for the standard nonlinear Schr\"odinger equation, solutions may disperse along non-free channels; see \cite{SW20223}. Building on the method initiated in \cite{SW20221}, this work characterizes the free part of the global solution, which we view as a first step toward establishing the soliton-resolution conjecture.

\subsection{Notations}\label{secF} Let $\bar{F}_c(\lambda)$, $F_j(\lambda)(j=1,2)$ denote smooth cut-off functions of the interval $[1,+\infty)$:
\begin{subequations}\label{F}
    
\begin{align}F_j(\lambda>a):=F_j(\lambda/a),\quad  \bar{F}_j(\lambda\leq a):=1-F_j(\lambda/a), \quad j=c,1,2,  
\end{align}
where $F_j(k), j=c,1,2,$ satisfy:
\be
F_j(k)=\begin{cases}
    1 & \text{ when }k\geq 1\\
    0 & \text{ when }k \leq 1/2
\end{cases}, \quad j=c,1,2.
\ee
Here it is worth noting that $F_j(k)$ is monotone increasing.

\end{subequations}

Throughout the paper, we use the notation $A \lesssim_s B$ and $ A \gtrsim_s B$ to indicate that there exists a
constant $C = C(s) > 0$ such that $A \leq CB$ and $A \geq CB,$ respectively.

\subsection{Adapted free channel wave operators}
Denote the evolution operator of the above system, from time $0$ to time $t$, as $U(t,0)$. The free channel wave operator associated with $H_0$ and $H:=H_0+V(x,t)$ is defined by 
\be
\Omega_\pm:=s\text{-}\lim\limits_{t\to \pm \infty} U(0,t)e^{-itH_0}\quad \text{ on }L^2_x(\mathbb{R}^n)
\ee
and its adjoint, by
\be
\Omega_\pm^*:=s\text{-}\lim\limits_{t\to \pm \infty} e^{itH_0}U(t,0)P_c\quad \text{ on }L^2_x(\mathbb{R}^n),\label{oWO}
\ee
where $P_c$ denotes the projection on the space of all free scattering states. Here, the term \emph{free-channel wave operator} is borrowed from multichannel scattering theory; likewise, our starting point—equation \eqref{Pc} below—also originates in that framework.
\begin{remark}
    When $V(x,t)$ is real and localized in $x$, the existence of $\Omega_\pm$ is classical; see Kato \cite{K1965} for Schr\"odinger equations. For bi-Laplacian Schr\"odinger equations, the existence of $\Omega_\pm$ follows by a similar argument (using unitarity of $U(0,t)$ and $e^{-itH_0}$ together with local decay estimates of the free flow). 
\end{remark}
\begin{remark}
The existence of $\Omega_\pm^{*}$ can be deduced from the definition of $P_c$. However, in time-dependent or nonlinear settings, the space of all free scattering states is not a-priori known. In this paper we construct the \emph{free-channel wave operator} following \cite{SW20221}, thereby bypassing any prior knowledge of $P_c$; once these wave operators are in hand, $P_c$ can then be defined from them provided that the system is linear and the Hamiltonian is self-adjoint.
\end{remark}

 In this paper, we stick to $\Omega_+^*$. We employ the technique developed in \cite{SW20221} to construct an adapted free channel wave operator. This operator coincides with $\Omega_\pm^*$, defined in \eqref{oWO}, in the weak sense. Indeed,
\be
w\text{-}\lim\limits_{t\to \infty}e^{itH_0}\psi(t)\quad \text{ in }L^2_x(\mathbb{R}^n)
\ee
exists whenever $\mathcal{N}(x,t,|\psi(t)|)$ and $\psi(t)$ satisfy our assumptions.

To be precise, we let $p:=-i\nabla_x$. The key observation is that 
\be
w\text{-}\lim\limits_{t\to \infty} e^{itH_0}(1-J(t,x,p))\psi(t)=0\quad \text{ in }L^2_x(\mathbb{R}^n),
\ee
provided that~\eqref{eq: uniformL2} holds. Then $\Omega_+^*\psi(0)$, defined in \eqref{oWO}, is equal to
\be
s\text{-}\lim\limits_{t\to +\infty} e^{itH_0}J(t,x,p)\psi(t),\quad \text{ in }L^2(\mathbb{R}^n),\label{1eq1}
\ee
and when $\mathcal N(x,t,|\psi(t)|)=V(x,t)$ is linear and real-valued such that $H(t)=H_0+V(x,t)$ is self-adjoint on the domain of $H(t)$, the projection on the space of all scattering states $P_c$ is given by 
\be
P_c=s\text{-}\lim\limits_{t\to\infty} U(0,t)J(t,x,p)U(t,0),\quad \text{ in }L^2_x(\mathbb{R}^n).\label{Pc}
\ee
In \cite{SW20221}, Soffer and Wu proposed using $J=F_c(\frac{|x-2tp|}{t^\alpha}\leq 1),  \alpha \in (0, 1-2/n)$ for cases where $V(x,t)\psi(t)\in L^\infty_tL^1_x(\R^{n+1})$, $n\geq 3$. For scenarios where $\langle x\rangle^\sigma V(x,t)\in L^\infty_{x,t}(\R^{n+1})$ for some $\sigma>1 $ and $ n\geq 1$, they recommended the use of $F_c( \frac{|x-2tp|}{t^\alpha}\leq 1)F_1(t^\beta|p|>1)$ for some $\alpha,\beta\in (0,1)$. 
 
\subsection{Main results} 
  In our setting, when $n\geq 1$, we assume that $ \langle x\rangle^\sigma V(x,t)\psi(t)\in L^\infty_t L^2_x(\mathbb{R}^n\times \mathbb{R}) $ for some $\sigma>1$. We take
\be
J(t,x,p)=e^{-itH_0}F_c(\frac{|x|}{t^\alpha}\leq 1)F_1(t^b|p|>1)e^{itH_0}
\ee
for all $b\in (0,(1/\sigma+1)/6)$ and $\alpha \in (b,(1+1/\sigma)/2)$. The adapted free channel wave operator acting on $\psi_0$, $\Omega_{\alpha,b}^*\psi_0$, is defined by
\be
\Omega_{\alpha,b}^{*}\psi_0=s\text{-}\lim\limits_{t\to \infty} \Omega_{\alpha,b}^*(t)\psi_0,\quad \text{ in }L^2_x(\mathbb{R}^n)\label{FNWO}
\ee
with 
\be
\begin{split}
 \Omega_{\alpha,b}^*(t)\psi_0:=&e^{itH_0}J(t,x,p)\psi(t)=F_c(\frac{|x|}{t^\alpha}\leq 1)F_1(t^b|p|>1)e^{itH_0}\psi(t).
 \end{split}
\ee

To state our main results, we also need the notion of weakly localized part of a solution,  introduced in~\cite{SW20221}.
\begin{definition} [The weakly localized part of the solution]~\cite[Definition 2.5]{SW20221}\label{def: weak} We say that a component  $\psi_{wloc}(t)$ of the solution is {\emph{weakly localized }} if it has non-zero mass and spreads slowly in the following sense: there exists $\beta\in (0,1)$ such that, for all $t\geq 1$,
\begin{equation}\label{eq: defweak}
    (\psi_{wloc}(t), |x|\psi_{wloc}(t))_{L^2_x(\mathbb{R}^n)}\lesssim t^\beta
\end{equation}
holds true.
\end{definition}

\begin{theorem}\label{Thm1}
Let $n\geq 1$, let $\sigma$ be as in Assumption \ref{condition}, and assume that  $V(x,t)$ and $\psi(t)$ satisfy Assumptions \ref{asp:sol} and \ref{condition}. Then:
 \begin{enumerate}
     \item  For all $b\in (0,(1/\sigma+1)/6)$ and $\alpha \in (b,(1+1/\sigma)/2)$, the adapted free channel wave operator acting on $\psi_0$, given by
\begin{align}\label{wlocalize1}
	\Omega_{\alpha,b}^*\psi_0:=s\text{-}\lim_{t\rightarrow\infty}\Omega_{\alpha,b}^*(t)\psi_0,
\end{align}
exists in $L_x^2(\mathbb{R}^n)$, and 
\eq
w\text{-}\lim\limits_{t\to \infty} \left(1-F_c(\frac{|x|}{t^\alpha}\leq 1)F_1(t^b|p|>1)\right)e^{itH_0}\psi(t)=0\label{wlimit}
\eeq
in $L_x^2(\mathbb{R}^n)$.
\item Moreover, if $\sigma
\geq 4$, the  following asymptotic decomposition holds: for all $\epsilon\in (0,1/4)$,
\be\label{wlocalize2}
\lim_{t\rightarrow\infty}\|\psi(t)-e^{-itH_0}\phi_+ -\psi_{wloc}(t)\|_{L^2_x(\mathbb{R}^n)}=0
\ee
where $\phi_+=\Omega^*_{\alpha,b}\psi_0\in L^2_x(\mathbb R^n)$ and 
\begin{equation}\label{eq: def weakloc}
    \psi_{wloc}(t)\equiv \psi_{wloc,\epsilon}(t):=F_c(\frac{|x|}{(t+1)^{\frac{1}{4}+\epsilon}}<1)\psi(t)
\end{equation}
is the weakly localized part of the solution, satisfying \eqref{eq: defweak} with   $\beta=\frac{1}{4}+\epsilon$.
\item In particular, when $n\geq 9$ and $\sigma> \frac{n}{2}$, the weakly localized part $\psi_{wloc}(t)$ defined in~\eqref{eq: def weakloc} is essentially localized in space in the following sense: there exists $\psi_{loc}(t)\in L^2_x(\mathbb R^n)$ (see~\eqref{eq: loc}) such that
\begin{equation}
    \|\psi_{wloc}(t)-\psi_{loc}(t)\|_{L^2_x(\mathbb R^n)}\to 0\qquad \text{ as }t\to \infty
\end{equation}
and
\begin{equation}
    \|\langle x\rangle^\delta \psi_{loc}(t)\|_{L^2_x(\mathbb R^n)}\leq C(\psi_0)
\end{equation}
for all $\delta\in [0,\frac{n}{8}-1)$ and some positive constant $C(\psi_0)>0$ depending on the initial datum $\psi_0$. 

\end{enumerate}
\end{theorem}

\begin{remark}
From~\eqref{wlimit}, we see that the wave operator $\Omega_+^*$ agrees with $\Omega_{\alpha,b}^*$.  
The cut-off $F_c\!\left(\tfrac{|x|}{t^\alpha} \leq 1\right)$, viewed as a propagation observable, confines the free evolution to a thin region of phase space and will be used to derive a propagation estimate.  
In low spatial dimensions, however, the time-dependent Schr\"odinger operator $H$ may have both a resonance and an eigenvalue at zero energy.  
To address this, we restrict the Schr\"odinger group $e^{itH_0}$ to the nonzero-frequency regime by inserting the cut-off $F_1(t^\beta |p| > 1)$.
\end{remark}
 
\begin{remark}Separating the free wave from the solution is the first step toward addressing the soliton resolution conjecture.  
In order to make progress on this conjecture, it is essential to understand certain properties of the ``non-free'' component.  
Theorem~\ref{Thm1} demonstrates how to separate the free wave by constructing the free channel wave operator, and further suggests that the weakly localized component (the ``non-free'' part) cannot spread faster than $t^{1/4+\epsilon}$, provided the interaction $V(x,t)$ is sufficiently localized in $x$.

\end{remark}

When $n\geq 5$, set
\be
J(t,x,p)=e^{-itH_0}F_c(\frac{|x|}{t^\alpha}\leq 1)e^{itH_0},\qquad  \alpha\in (0,\frac{1}{2}-\frac{2}{n})
\ee
and define 
\begin{align}
    \Omega_{\alpha}^{*}(t)\psi(0):=F_c(\frac{|x|}{t^\alpha}\leq 1)e^{itH_0}\psi(t).
\end{align}
Then we have following result:
\begin{theorem}\label{Thm2}
If Assumptions \ref{asp:sol} and~\ref{condition1} hold, then for $n\geq 5$ and every $\alpha\in (0,\frac{1}{2}-\frac{2}{n}),$ the adapted channel wave operator acting on $\psi_0$, defined by
\begin{align}\label{exist-Ome}
	\Omega_{\alpha}^{*}\psi_0:=s\text{-}\lim_{t\rightarrow\infty}\Omega_{\alpha}^{*}(t)\psi_0,\ 
\end{align}
exists in  $L_x^2(\mathbb{R}^n).$ Furthermore, $\Omega_\alpha^*\psi_0$ is independent on the choice of $\alpha$ in the following sense: for all $\alpha,\alpha'\in (0, \frac{1}{2}-\frac{2}{n})$, 
\eq\label{chooseindepen}
\Omega_{\alpha}^*\psi_0=\Omega_{\alpha'}^*\psi_0.
\eeq
\end{theorem}
\begin{remark} It is worth noting that when $n \geq 5$, the terms $\mathcal{N}(x,t,|\psi(t)|)\psi(t)$ or $V(x,t)\psi(t)$ may be nonlocal in the spatial variable $x$.  
The condition $n \geq 5$ corresponds to the requirement that 
\[
\|e^{itH_0}\|_{L^1_x(\mathbb{R}^n)\to L^\infty_x(\mathbb{R}^n)} \in L^1_t[1,\infty),
\]
as stated in~\eqref{decayestimate}.

\end{remark}

For notational convenience, we suppress the arguments of $F_c$ and $F_1,$  writing $F_c\equiv F_c(\frac{|x|}{t^\alpha}\leq 1)$ and  $F_1\equiv F_1(t^b|p|>1)$ when  the meaning is clear.

\section{Estimates for interaction terms and Commutator Estimates}
In this section, we collect key estimates for the free evolution, several results from~\cite{SW20221} and bounds for the interaction terms that will be used in the proofs of Theorems~\ref{Thm1} and~\ref{Thm2}.

\subsection{Decay estimates of the free flow}
For free bi-Laplacian Schr\"odinger equations, the  propagator $e^{-it\Delta^2}$  can be written as 
\begin{align}
    e^{-it\Delta^2}f=\F^{-1}(e^{-it|q|^4}\hat{f})=\int_{\mathbb{R}^n}I_0(t,x-y)f(y)dy,
\end{align}
where $I_0(t,x)=\F^{-1}(e^{-it|q|^4})(x)$ is the distributional kernel of $e^{-it\Delta^2}.$
For biharmonic operator $(-\Delta)^2$, Ben-Artzi, Koch and Saut \cite{MR1745182}
proved the following sharp kernel estimate
\begin{align}
    |D^{\gamma}I_0(t,x)|\leq C|t|^{-(n+|\gamma|)/4}(1+|t|^{-1/4}|x|)^{(|\gamma|-n)/3},\quad t\not=0, x\in \mathbb{R}^n,
\end{align}
which implies
\begin{align}\label{decay1}
    \|D^{\gamma}e^{-itH_0}\|_{L^1_x(\mathbb R^n)\rightarrow L^{\infty}_x(\mathbb R^n)}\lesssim|t|^{-\frac{n+|\gamma|}{4}}, \quad |\gamma|\leq n.
\end{align}
Here $D=(\partial_{x_1},\partial_{x_2},\cdots,\partial_{x_n})$, $D^\gamma=\partial_{x_1}^{\gamma_1}\cdots \partial_{x_n}^{\gamma_n}$ for $\gamma=(\gamma_1,\cdots,\gamma_n)\in \mathbb N^n$ and $|\gamma|=\sum\limits_{j=1}^n\gamma_j$. Specializing \eqref{decay1} to $\gamma=0$ yields the   $L^1_x\rightarrow L^{\infty}_x$ dispersive estimate for    free bi-Laplacian   Schr\"odinger equation:
\begin{align}\label{decayestimate}
     \|e^{-itH_0}\|_{L^1_x(\mathbb R^n)\rightarrow L^{\infty}_x(\mathbb R^n)}\lesssim|t|^{-\frac{n}{4}}.
\end{align}
\subsection{Propagation estimates}\label{sec: PP}
Throughout the paper, we   use the following tools introduced in \cite{SW20221}.
\begin{enumerate}
 \item(\textbf{Propagation Estimate}) Given a bounded operator $B$, we define
\begin{align}
    \langle B(t): \psi(t)\rangle_t:=(\psi(t), B(t)\psi(t))_{L^2_x(\mathbb{R}^n)}=\int_{\mathbb{R}^3}\left(\psi(t)\right)^* B(t)\psi(t)dx,
\end{align}
where $\psi(t)$ denotes the solution to \eqref{fourth-order} and $f^*$ stands for the complex conjugation of  $f$. We call the family  $\{B(t)\}_{t\geq1}$  a \textbf{Propagation Observable  } if  
\begin{align}
    & \partial_t\langle B(t):\psi(t)\rangle_t=(\psi(t),C^*C\psi(t))_{L^2_x(\mathbb R^n)}+g(t),\nonumber\\
    & g(t)\in L^1_{t}[1,\infty),\quad C^*C\geq0.
\end{align}
Assume  Condition~\ref{eq: uniformL2} holds. Then, for any bounded operator and $T>t_0\geq 1$,
\begin{equation}
    \sup\limits_{t\in [t_0,T]} \left|( \psi(t), B(t)\psi(t))_{L^2_x(\mathbb{R}^n)}\right|<\infty.
\end{equation}
Integrating in time yields the following \textbf{Propagation Estimate}: for $T>t_0\geq 1$,
\begin{align}
    \int_{t_0}^T\|C(t)\phi(t) \|_{L^2_x(\mathbb{R}^n)}^2dt&=( \psi(T), B(T)\psi(T))_{L^2_x(\mathbb{R}^n)}-\nonumber\\
    &( \psi(t_0), B(t_0)\psi(t_0))_{L^2_x(\mathbb{R}^n)}-\int_{t_0}^Tg(s) ds\nonumber\\
&\quad\quad\leq 2\sup\limits_{t\in [t_0,T]} \left|( \psi(t), B(t)\psi(t))_{L^2_x(\mathbb{R}^n)}\right|+C_g,  
\end{align}
where $C_g:=\|g(t)\|_{L^1_t[1,\infty)}.$ 
 \item(\textbf{Relative Propagation Estimate})  Given a family of bounded operators $\{\tilde{B}(t)\}_{t\geq 1}$,  define 
 \begin{align}
     \langle \tilde{B}: \phi(t)\rangle_t:=(\phi(t), \tilde{B}(t)\phi(t)  )_{L^2_x(\mathbb{R}^n)}=\int_{\mathbb{R}^n}\left(\phi(t)\right)^*\tilde{B}(t)\phi(t)dx,
 \end{align}
 where $\phi(t)$ is a flow related to $\psi(t)$ such that
\eq
\sup\limits_{t\geq 1}\langle \tilde{B}:\phi(t)\rangle_t<\infty. \label{phiH}
\eeq
Assume further that 
\begin{align}\label{eq: rpp C1}
&\partial_t\langle \tilde{B} : \phi(t)\rangle_t=\pm\sum\limits_{l=1}^L\langle \phi(t), C_l^*C_l\phi(t)\rangle+g(t),\nonumber\\
&g(t)\in L^1[1,\infty), \quad C_l^*C_l\geq 0,\quad l=1,\cdots,L .
\end{align}
We then refer to  $\{\tilde{B}(t)\}_{t\geq 1}$ as a {\bf Relative Propagation Observable} with respect to $\phi(t)$.
Integrating in time gives, for $l=1,\cdots,L$ and $1\leq t_0<T$,
\eq
\int_{t_0}^T\|C_l(t)\phi(t) \|_{L^2_x(\mathbb{R}^n)}^2dt\leq 2\sup\limits_{t\in [t_0,T]} \left|( \phi(t), \tilde{B}(t)\phi(t))_{L^2_x(\mathbb{R}^n)}\right|+C_g, \label{CC0}
\eeq
where $ C_g:=\|g(t)\|_{L^1_t[1,\infty)}.$
We refer to  \eqref{CC0} as a {\bf Relative Propagation Estimate}.
\end{enumerate}
\subsection{Commutator estimates}
We include here the commutator estimate~\cite{SW20221}.
\begin{lemma}\label{commutatorlemma}
	For $t\geq1,$ $b<\alpha \leq 1$, $l=0,1$  we have 
	\begin{align}\label{commutator1}
		\|[F_c,F_1^{(l)}]\|_{L^2_x(\mathbb{R}^n)\rightarrow L^2_x(\mathbb{R}^n)}\lesssim\frac{1}{t^{\alpha-b}},
	\end{align}
\begin{align}\label{commutator2}
		\|[F^{(l)}_c,F_1]\|_{L^2_x(\mathbb{R}^n)\rightarrow L^2_x(\mathbb{R}^n)}\lesssim\frac{1}{t^{\alpha-b}},
\end{align}
where 
\begin{align}
	F^{(l)}_1(k):=\frac{d^l}{dk^l}[F_1] \quad \mbox{and}\quad 	F^{(l)}_c(k):=\frac{d^l}{dk^l}[F_c].
\end{align}
\end{lemma}
Write 
\begin{equation}
    \partial_t{F_1(t^b |p|>1)}=\frac{1}{t} \left(t \partial_t{F_1(t^b |p|>1)}\right).
\end{equation}
We view $t \partial_t{F_1(t^b |p|>1)}$ as an operator similar to $F_1(t^b|p|>1)$ in Lemma~\ref{commutatorlemma} and then use Lemma~\ref{commutatorlemma} to obtain
\begin{corollary}\label{Cor1}
    For $t\geq1,$ $b<\alpha \leq 1$,   we have 
\begin{align}
    \|[F_c,\partial_t[F_1]]\|_{L^2_x(\mathbb{R}^n)\rightarrow L^2_x(\mathbb{R}^n)}\lesssim\frac{1}{t^{\alpha-b+1}}
\end{align}
and 
\begin{align}
    \|[\partial_t[F_c],F_1]\|_{L^2_x(\mathbb{R}^n)\rightarrow L^2_x(\mathbb{R}^n)}\lesssim\frac{1}{t^{\alpha-b+1}}.
\end{align}
\end{corollary}
\subsection{Estimates for the interaction terms} We need the following estimates for the interaction terms.

\begin{proposition}\label{prop2.1}
Let $n\geq 1$, and assume $V(x,t)$ and $\psi(t)$ satisfy Assumption \ref{condition}, with $\sigma$   as therein. For  any $b\in (0,(1/\sigma+1)/6)$ and $\alpha \in (0,(1+1/\sigma)/2)$,  we have, for all $t\geq1,$
    \begin{align} \label{interaction1}
    \|F_cF_1e^{itH_0}V(x,t)\psi(t)\|_{L^2_x(\mathbb{R}^n)}\lesssim_{\sigma} \frac{1}{t^{(\sigma+1)/2}}\|V(x,t)\psi(t)\|_{L^{\infty}_t L^2_{x,\sigma}(\mathbb{R}^{n+1})}.
\end{align}
\end{proposition}
\begin{proof}
Let
    \begin{align}
	a_{\psi}(t):=(-i)F_cF_1e^{itH_0}V(x,t)\psi(t)
\end{align}
and set  $e=\frac{1}{2}(1+\frac{1}{\sigma})$. We split $a_\psi(t)$ as
\begin{align}
    a_{\psi}(t)&=(-i)F_cF_1e^{itH_0}\chi(|x|>t^e)V(x,t)\psi(t)\nonumber\\
    &+(-i)F_cF_1e^{itH_0}\chi(|x|\leq t^e)V(x,t)\psi(t)=:a_{\psi,1}(t)+a_{\psi,2}(t).
\end{align}
Using Assumption~\ref{condition},  we obtain the time decay for $a_{\psi,1}(t)$:
\begin{align}\label{a1}
    \|a_{\psi,1}(t)\|_{L^2_x}&\lesssim \frac{1}{t^{e\sigma}}\|V(x,t)\psi(t)\|_{L^{\infty}_t L^2_{x,\sigma}(\mathbb{R}^{n+1})}
\end{align}
with $e\sigma=(1+\sigma)/2$.\par
For $a_{\psi,2}(t),$ we apply the method of nonstationary phase (integration by parts in Fourier space, see also (B.5)–(B.16) in~\cite{SW20221}) to obtain
\begin{align}\label{a2}
    \|a_{\psi,2}(t)\|_{L^2_x}\lesssim_N\frac{1}{t^N}\|V(x,t)\psi(t)\|_{L^{\infty}_t L^2_{x,\sigma}(\mathbb{R}^{n+1})}.
\end{align}
More precisely, we use
\begin{align}
	F_c(\frac{|x|}{t^{\alpha}}\leq 1)&e^{i(x-y)\cdot q}e^{itq^4}\chi(|y|\leq t^e)=\nonumber\\
	&\frac{1}{i[(x-y)\cdot \hat{q}+4t|q|^3]}\partial_{|q|}[	F_c(\frac{|x|}{t^{\alpha}}\leq 1)e^{i(x-y)\cdot q}e^{itq^4}\chi(|y|\leq t^e)]
\end{align}
with
\begin{align}
	|(x-y)\cdot \hat{q}+4t|q|^3|\gtrsim t^{1-3b}\gtrsim t^e,\qquad t\geq 1,
\end{align}
when $b<(1-e)/3=(1/\sigma+1)/6$ and $\alpha<e=(1+1/\sigma)/2$. Choosing $N=[ \sigma]+1$, we obtain
\begin{align}
		\|a_{\psi,2}(t)\|_{L^2_x(\mathbb{R}^n)}\lesssim_{\sigma} \frac{1}{t^{e\sigma}}\|V(x,t)\psi(t)\|_{L^{\infty}_t L^2_{x,\sigma}(\mathbb{R}^{n+1})},
\end{align}
which, together with~\eqref{a1}, yields~\eqref{interaction1}. This completes the proof.
\end{proof}

\begin{proposition}\label{prop2.2}
Let $n\geq 1$. Assume Assumptions~\ref{asp:sol} and~\ref{condition} hold, and let $m$ and $\sigma$ be as in Assumptions~\ref{asp:sol} and \ref{condition}, respectively. For any $b\in (0,(1/\sigma+1)/6)$ and $\alpha \in (0,(1+1/\sigma)/2)$,  we have, for all $t\geq1,$
\begin{align}\label{interaction2}
	&	|(F_1e^{itH_0}\psi(t),F_cF_1e^{itH_0}V(x,t)\psi(t))_{L^2_x(\mathbb{R}^n)}|\nonumber\\
 \lesssim_{\sigma} &\frac{m}{t^{(\sigma+1)/2}}\|V(x,t)\psi(t)\|_{L^{\infty}_t L^2_{x,\sigma}(\mathbb{R}^{n+1})}
\end{align}
and
\begin{align}\label{interaction3}
		&|(F_1e^{itH_0}V(x,t)\psi(t),F_cF_1F_c e^{itH_0}\psi(t))_{L^2_x(\mathbb{R}^n)}|\nonumber\\
  \lesssim_{\sigma} &\frac{m}{t^{(\sigma+1)/2}}\|V(x,t)\psi(t)\|_{L^{\infty}_t L^2_{x,\sigma}(\mathbb{R}^{n+1})}.
\end{align}
\end{proposition}

\begin{proof}
Estimates~\eqref{interaction2} and~\eqref{interaction3} follow directly from Assumptions~\ref{asp:sol} and~\ref{condition}, together with Proposition~\ref{prop2.1}.
\end{proof}

\begin{proposition}\label{prop2.3}
Let $n\geq 5$,  and assume Assumption \ref{condition1} holds. Then for all $\alpha\in (0,\frac{1}{2}-\frac{2}{n})$ and $t\geq 1$, we have 
\begin{align}\label{interaction4}
    \|F_ce^{itH_0}V(x,t)\psi(t)\|_{L^2_x (\mathbb{R}^n)}\lesssim \frac{1}{t^{\frac{n}{4}(1-2\alpha)}}\|V(x,t)\psi(t)\|_{L^\infty_tL^1_x(\mathbb{R}^{n+1}) }.
\end{align}
\end{proposition}
\begin{proof}
   Combining the $L^\infty$ decay estimate for the free flow \eqref{decayestimate} with Assumption~\ref{condition1} to obtain
\begin{align}
    \|F_ce^{itH_0}V(x,t)\psi(t)\|_{L^2_x (\mathbb{R}^n)} \lesssim &\|F_c\|_{L^2_x(\mathbb{R}^n)}\|e^{itH_0}\|_{L^1_x\rightarrow L^{\infty}_x}\|V(x,t)\psi(t)\|_{L^1}\nonumber\\
    \lesssim &\frac{1}{t^{\frac{n}{4}(1-2\alpha)}}\|V(x,t)\psi(t)\|_{L^\infty_tL^1_x(\mathbb{R}^{n+1}) }.
\end{align}\end{proof}


\subsection{Idea of establishing the relative propagation estimates}\label{sec: idea PP} 
In order to establish the relative propagation estimates, we use the following identity to produce the positive term: for a bounded, positive operator $X$ and a bounded operator $Y$,
\begin{equation}\label{eq: XY}
XY+YX=2\sqrt{X}\,Y\,\sqrt{X}+[\sqrt{X},[\sqrt{X}, Y]],
\end{equation}
where $[X,Y]=XY-YX$ denotes the commutator. When $Y$ is also positive, the first term on the r.h.s. of \eqref{eq: XY} yields the positive operator appearing in the relative propagation estimates (see \eqref{eq: rpp C1} and \eqref{CC0}). For the remainder term (the second term on the r.h.s. of \eqref{eq: XY}), its expectation can typically be controlled via commutator bounds.

To be precise, in our setting, we take $F_c\equiv F_c(\frac{|x|}{t^\alpha}\leq 1)$ and $F_1\equiv F_1(t^b|p|>1)$, and observe 
\begin{equation}
    \langle B(t): \phi(t)\rangle_t,\qquad B(t):=F_1F_cF_1,\quad \phi(t)=e^{itH_0}\psi(t)
\end{equation}
for $t\geq 1$. To establish the propagation estimates, we compute 
\begin{equation}\label{eq: eq1}
\begin{aligned}
    &\partial_t\langle B(t): \phi(t)\rangle_t=\langle F_1\partial_t[F_c]F_1: \phi(t)\rangle_t+\langle \partial_t[F_1]F_cF_1+ F_1F_c\partial_t[F_1]:\phi(t)\rangle_t\\
    &+( \phi(t), (-i)F_1F_cF_1e^{itH_0}V\psi(t))_{L^2_x(\mathbb R^n)}+( (-i)F_1F_cF_1e^{itH_0}V\psi(t),\phi(t))_{L^2_x(\mathbb R^n)},
    \end{aligned}
\end{equation}
where the last two terms on the r.h.s. of \eqref{eq: eq1} are interaction terms, which can be controlled directly. Since $F_1\partial_t[F_c]F_1$ is bounded and positive, it suffices to rearrange the second term of the the r.h.s. of~\eqref{eq: eq1}. We find 
\begin{equation}\label{eq: eq2}
    \partial_t[F_1]F_cF_1+F_1F_c\partial_t[F_1]=\partial_t[F_1]F_1 F_c+F_cF_1\partial_t[F_1]+\left( \partial_{t}[F_1][F_c,F_1]+[F_1,F_c]\partial_t[F_1]\right).
\end{equation}
For the first two terms of the r.h.s. of~\eqref{eq: eq2}, by taking $Y=\partial_t[F_1]F_1$ and $X=F_c$ and by~\eqref{eq: XY}, we obtain 
\begin{equation}\label{eq: eq3}
\partial_t[F_1]F_1F_c+F_c\partial_t[F_1]F_1=2\sqrt{F_c} \partial_t[F_1]F_1\sqrt{F_c}+[\sqrt{F_c},[\sqrt{F_c},\partial_t[F_1]F_1]].
\end{equation}
By the Commutator estimates (Lemma~\ref{commutatorlemma} and Corollary~\ref{Cor1}), the operator norms of the commutator terms appearing in \eqref{eq: eq2} and \eqref{eq: eq3} are absolutely integrable in time. We refer to these as remainder operators, and their expectations as remainder terms.
\subsection{Estimates for the reminder terms} Building on Section~\ref{sec: idea PP}, we now estimate the remainder terms needed to establish the relative propagation estimates.  We list in this section the reminder terms we need and their estimates. 

Fix $\phi(t)=e^{itH_0}\psi(t)$. To observe $\langle B_1(t): \phi(t)\rangle_t$ with $B_1(t)=F_1F_cF_1$, we have the following reminder term:
\begin{align}\label{r}
    \langle r(t):\phi(t)\rangle_t:&=\left(\phi(t),(\partial_{t}[F_1][F_c,F_1]+[F_1,F_c]\partial_t[F_1])\phi(t)\right)_{L^2_x(\R^n)}   \nonumber\\
    &+ \left(\phi(t),[\sqrt{F_c},[\sqrt{F_c},\partial_t[F_1]F_1]]\phi(t)\right)_{L^2_x(\R^n)}.
\end{align} 
To observe $\langle B_2(t): \phi(t)\rangle_t$ with $B_1(t)=F_cF_1F_c$, we have the following reminder term:
\begin{align}\label{r-tilde}
   \langle \tilde r(t):\phi(t)\rangle_t:&=\left(\phi(t),(\partial_{t}[F_c][F_1,F_c]+[F_c,F_1]\partial_t[F_c])\phi(t)\right)_{L^2_x(\R^n)}   \nonumber\\
    &+ \left(\phi(t),[\sqrt{F_1},[\sqrt{F_1},\partial_t[F_c]F_c]]\phi(t)\right)_{L^2_x(\R^n)}.
\end{align} 

\begin{lemma}\label{lem2.5}
   Assume Assumption~\ref{asp:sol} holds and let $m$ be as in~\eqref{eq: uniformL2}. Then, for all $t\geq1$, $0<b<\alpha$, we have
\begin{align}
    | \langle r(t):\phi(t)\rangle_t|\lesssim\frac{m^2}{t^{1+\alpha-b}}\label{com: eq1}
\end{align}
and 
\begin{align}
    | \langle \tilde r(t):\phi(t)\rangle_t|\lesssim\frac{m^2}{t^{1+\alpha-b}}.\label{com: eq2}
\end{align}
\end{lemma}
\begin{proof}\eqref{com: eq1} and \eqref{com: eq2} follow directly from the definitions of $r(t)$ and $\tilde r(t)$, together with Lemma~\ref{commutatorlemma} and Corollary~\ref{Cor1}.\end{proof}

\section{Proof of the existence of the adapted wave operators}
\subsection{Proof of Theorem~\ref{Thm2}}
We now use the \textbf{Relative Propagation Estimate} introduced in Section \ref{sec: PP} to prove Theorem~\ref{Thm2}. In this setting we fix the spatial dimension $n\ge 5$.

\begin{proof}[{Proof of Theorem \ref{Thm2}.}]
For convenience, we adopt the convention $F_c\equiv F_c(\frac{|x|}{t^\alpha}\leq 1)$ or $F_c\equiv F_c(\frac{|x|}{s^\alpha}\leq 1)$ when the context is clear. We define 
\begin{align}
   \Omega_\alpha^*(t)\psi_0=F_ce^{itH_0}\psi(t)
\end{align}
and by Assumption~\ref{asp:sol}, $\Omega_\alpha^*(t)\psi_0\in L^2_x(\mathbb{R}^n)$ for all $t\geq 1$. We apply
Cook's method to expand $ \Omega_{\alpha}^*(t)\psi_0$ for $t\geq1$:
\begin{align}
\Omega_\alpha^*(t)\psi_0=&\Omega_\alpha^*(1)\psi_0+(-i)\int_1^tds F_ce^{isH_0}V(x,s)\psi(s)+\int_1^tds    \partial_s[F_c]e^{isH_0}\psi(s)\nonumber\\
=:&\Omega_\alpha^*(1)\psi_0+\Omega_{\alpha,1}(t)\psi_0+\Omega_{\alpha,2}(t)\psi_0.
\end{align}
Thanks to Assumption \ref{asp:sol} and the unitarity of $e^{itH_0}$, 
\eq
\|\Omega_\alpha^*(1)\psi_0\|_{L^2_x(\mathbb{R}^n)}\leq\sup\limits_{t\in\R} \|\psi(t)\|_{L^2_x(\mathbb{R}^n)}=m.\label{bounded1}
\eeq
For the interaction term $\Omega_{\alpha,1}(t)$, we use Proposition \ref{prop2.3} to obtain the existence of
  \begin{align}\label{bounded3}
     \Omega_{\alpha,1}\psi_0:=\lim_{t\rightarrow\infty}\Omega_{\alpha,1}(t)\psi_0
 \end{align}
in $L^2_x(\mathbb{R}^n).$ 

For $\Omega_{\alpha,2}(t)\psi_0$, we use the \textbf{Relative Propagation Estimate} to prove $\{\Omega_{\alpha,2}(t)\psi_0\}_{t\geq 1}$ is Cauchy in $L^2_x(\mathbb R^n)$. By taking
$$
\begin{cases}
B_1(t)=F_c,\\
\phi(t)=e^{itH_0}\psi(t),
\end{cases}
$$ 
we obtain
\eq
\langle B_1(t): \phi(t)\rangle_t\leq \left(\sup\limits_{t\in \R}\|\psi(t)\|_{L^2_x(\mathbb{R}^n)}\right)^2=m^2<\infty.
\eeq
Using
\begin{align}
\partial_tF_c(\frac{|x|}{t^\alpha}\leq 1)=\frac{-\alpha}{t}\times \left(F_c'(\frac{|x|}{t^\alpha}\leq 1)\frac{|x|}{t^{\alpha}}\right)\geq0,
\end{align}
we compute
\begin{align}
\partial_t\langle B(t):\phi(t)\rangle_t=b_1(t)+b_2(t)+b_3(t),
\end{align}
where $b_j(t), j=1,2,3$ are given by
\begin{align}
b_1(t):=(e^{itH_0}\psi(t),\partial_tF_c(\frac{|x|}{t^\alpha}\leq 1)e^{itH_0}\psi(t))\geq 0,
\end{align}
\begin{align}
b_2(t):=(-i)(e^{itH_0}\psi(t), F_c(\frac{|x|}{t^\alpha}\leq 1)e^{itH_0}V(x,t)\psi(t))
\end{align}
and
\begin{align}
b_3(t):=i(e^{itH_0}V(x,t)\psi(t), F_c(\frac{|x|}{t^\alpha}\leq 1)e^{itH_0}\psi(t)),
\end{align}
respectively. By H\"{o}lder's inequality, Assumption~\ref{asp:sol} and Proposition \ref{prop2.3}, we obtain 
 \begin{align}
     |b_2(t)+b_{3}(t)|\lesssim\frac{m}{t^{\frac{n}{4}(1-2\alpha)}}\|V(x,t)\psi(t)\|_{L^\infty_tL^1_x(\mathbb{R}^n) }\in L^1_t[1,\infty),
 \end{align}
provided that $\alpha\in (0,\frac{1}{2}-\frac{2}{n})$ and $n\geq5$. Then by the \textbf{Relative Propagation Estimate} \eqref{CC0}, we have, for all $T\geq 1$, 
 \begin{align}
     \int_1^T b_1(t)dt &\leq \langle B(t): \phi(t)\rangle_t|_{t=T}-\langle B(t): \phi(t)\rangle_t|_{t=1}+\|b_2(t)+b_3(t)\|_{L^1_t[1,\infty)}\nonumber\\
     &\leq 2m^2+\|b_2(t)+b_3(t)\|_{L^1_t([1,\infty))}.
 \end{align}
 Hence, 
 \begin{align}
  \lim\limits_{T\to \infty}\int_1^{T} b_1(t)dt  \text{ exists and } \int_1^{\infty} b_1(t)dt<\infty.
 \end{align}
 By H\"{o}lder's inequality in $s$ variable, we have, for $t_2\geq t_1\geq T\geq1$,
 \begin{align*}
     \|\Omega_{\alpha,2}(t_2)\psi_0&-\Omega_{\alpha,2}(t_1)\psi_0\|_{L^2_x(\mathbb{R}^n)}=\left\|\int_{t_1}^{t_2}    \partial_s[F_c]e^{isH_0}\psi(s)ds \right\|_{L^2_x(\mathbb{R}^n)}\\
     &\leq \left\|\left(\int_{t_1}^{t_2}    \partial_s[F_c]ds\right)^{1/2}\left(\int_{t_1}^{t_2}    \partial_s[F_c]|e^{isH_0}\psi(s)|^2ds\right)^{1/2}\right\|_{L^2_x(\mathbb{R}^n)}\\
     &\leq \left(\int_T^{\infty}b_1(s)ds\right)^{1/2}\rightarrow0,
 \end{align*}
 as $T\rightarrow\infty.$ Hence, $\{\Omega_{\alpha,2}(t)\}_{t\geq 1} $ is Cauchy   in $L^2_x(\mathbb{R}^n)$ and subsequently 
 \begin{align}\label{bounded2}
     \Omega_{\alpha,2}\psi_0:=\lim_{t\rightarrow\infty}\Omega_{\alpha,2}(t)\psi_0\text{
 exists in $L^2_x(\mathbb{R}^n).$}
 \end{align}
Combining \eqref{bounded1},\eqref{bounded3} with \eqref{bounded2}, yields \eqref{exist-Ome}.
 
 We also have that for all $\alpha,\alpha'\in (0,\frac{1}{2}-\frac{2}{n})$ and $\phi\in L^2_x(\mathbb{R}^n)$, by Cauchy-Schwarz inequality, Assumption~\ref{asp:sol} and the unitarity of $e^{itH_0}$,
\begin{align}
    \left|(\phi, \Omega_\alpha^*(t)\psi_0-\Omega_{\alpha'}^*(t)\psi_0)_{L^2_x(\mathbb{R}^n)}\right|=& \left|( (F_c(\frac{|x|}{t^\alpha}\leq 1)-F_c(\frac{|x|}{t^{\alpha'}}\leq 1))\phi, e^{itH_0}\psi(t))_{L^2_x(\mathbb{R}^n)}\right|\nonumber\\
    \leq & \|(F_c(\frac{|x|}{t^\alpha}\leq 1)-F_c(\frac{|x|}{t^{\alpha'}}\leq 1))\phi \|_{L^2_x(\mathbb{R}^n)}\| e^{itH_0}\psi(t)\|_{L^2_x(\mathbb{R}^n)}\nonumber\\
    \leq & m\|(F_c(\frac{|x|}{t^\alpha}\leq 1)-F_c(\frac{|x|}{t^{\alpha'}}\leq 1))\phi \|_{L^2_x(\mathbb{R}^n)}\nonumber\\
    \to & 0
\end{align}
as $t\to \infty$. This implies 
\eq
w\text{-}\lim\limits_{t\to \infty} \Omega_\alpha^*(t)\psi_0-\Omega_{\alpha'}^*(t)\psi_0=0,\quad \text{ in }L^2_x(\mathbb{R}^n)
\eeq
and therefore, due to the existence of $\Omega_\alpha^*\psi_0$ and $\Omega_{\alpha'}^*\psi_0$ in $L^2_x(\mathbb{R}^n)$ in strong sense, Eq.~\eqref{chooseindepen} holds. We conclude the proof of Theorem \ref{Thm2}.

\end{proof}
\subsection{Proof of~Theorem~\ref{Thm1} part \emph{(i)}}

For lower spatial dimensions ($n\ge 1$), a similar argument applies.

 \begin{proof}[{Proof of Theorem \ref{Thm1} part (i).}]
  For spatial dimension $n\ge 1$, we define
 \begin{align}
     \Omega_{\alpha,b}^*(t)\psi_0:= F_c F_1 e^{itH_0}\psi(t)
 \end{align}
 and by Assumption~\ref{asp:sol}, $\Omega_{\alpha,b}^*(t)\psi_0\in L^2_x(\mathbb{R}^n)$ for all $t\geq 1$. We now apply Cook’s method to expand $\Omega_{\alpha,b}^*(t)\psi_0$  for $t\geq1$:
\begin{align}
	\Omega_{\alpha,b}^*(t)\psi_0&=\Omega_{\alpha,b}^*(1)\psi_0+(-i)\int_{1}^t F_c F_1 e^{isH_0}V(x,s)\psi(s)\,ds \nonumber\\
	&\quad +\int_{1}^t\partial_s\!\left[F_c \right]F_1 e^{isH_0}\psi(s)\,ds+\int_{1}^t F_c \partial_s\!\left[F_1 \right]e^{isH_0}\psi(s)\,ds  \nonumber\\
	& =:\Omega_{\alpha,b}^*(1)\psi_0+\psi_{\Omega,1}(t)+\psi_{\Omega,2}(t)+\psi_{\Omega,3}(t).
\end{align}
Thanks to Assumption \ref{asp:sol} and the unitarity of $e^{itH_0}$, we also have 
\eq
\|\Omega_{\alpha,b}^*(1)\psi_0\|_{L^2_x(\mathbb{R}^n)}\leq\sup\limits_{t\in\R} \|\psi(t)\|_{L^2_x(\mathbb{R}^n)}=m. 
\eeq

For $\psi_{\Omega,1}(t)$,
 by Proposition \ref{prop2.1} and the condition that $\sigma>1$, we have 
\begin{equation}     \|F_cF_1e^{isH_0}V(x,s)\psi(s)\|_{L^2_x(\mathbb R^n)}\in L^1_s[1,\infty)
 \end{equation} 
 and subsequently
\begin{align}\label{limit4}
		\psi_{\Omega,1}(\infty):=\lim_{t\rightarrow\infty}\psi_{\Omega,1}(t)\ \mbox{exists in }\ L_x^2(\mathbb{R}^n).
\end{align}
For $\psi_{\Omega,2}(t)$, we use the \textbf{Relative Propagation Estimate} to prove $\{\psi_{\Omega,2}(t)\}_{t\geq 1}$ is Cauchy in $L^2_x(\mathbb R^n)$. By taking 
$$
\begin{cases}
B(t)=F_1F_cF_1,\\
\phi(t)=e^{itH_0}\psi(t),
\end{cases}
$$ 
we obtain
\eq
\langle B(t): \phi(t)\rangle_t\leq \left(\sup\limits_{t\in \R}\|\psi(t)\|_{L^2_x(\mathbb{R}^n)}\right)^2=m^2<\infty.
\eeq
Using
\begin{align}
\partial_tF_c(\frac{|x|}{t^\alpha}\leq 1)=\frac{-\alpha}{t}\times \left(F_c'(\frac{|x|}{t^\alpha}\leq 1)\frac{|x|}{t^{\alpha}}\right)\geq0,
\end{align}
we compute   
\eq
\partial_t[\langle B(t):\phi(t)\rangle_t]=c_{1}(t)+c_{2}(t)+ \langle r(t):\phi(t)\rangle_t,
\eeq
where
\begin{align*}
c_{1}(t):=&(F_1 \phi(t), \partial_t[F_c]F_1 \phi(t))_{L^2_x(\R^n)}\geq 0,\quad \forall t\geq 1,
\end{align*}
\begin{align*}
   c_{2}(t):&= (-i)(F_1e^{itH_0}\psi(t), F_cF_1e^{itH_0} V(x,t)\psi(t))_{L^2_x(\R^n)}\\
&\quad+i(F_1e^{itH_0}V(x,t)\psi(t), F_cF_1e^{itH_0}\psi(t))_{L^2_x(\R^n)}
\end{align*}
 and  $\langle r(t):\phi(t)\rangle_t$  is defined in \eqref{r}.
 Then using  H\"{o}lder's inequality,  Proposition \ref{prop2.2} and Lemma \ref{lem2.5}, we get
 \begin{align}
    |c_2(t)|\lesssim \frac{m}{t^{(\sigma+1)/2}}\|V(x,t)\psi(t)\|_{L^{\infty}_t L^2_{x,\sigma}(\mathbb{R}^{n+1})}\in L^1_t[1,\infty)
 \end{align}
 and 
 \begin{align}
     |\langle r(t):\phi(t)\rangle_t|\lesssim\frac{1}{t^{1+\alpha-b}}m^2\in L^1_t[1,\infty),
 \end{align}
 provided that $\sigma>1$ and $\alpha>b$.
 
Then by the \textbf{Relative Propagation Estimate} \eqref{CC0}, we have, for all $T\geq 1$,
\begin{align}
\int_1^T c_{1}(t)dt
\leq &\langle B(t):\phi(t)\rangle_t\vert_{t=T}-\langle B(t):\phi(t)\rangle_t\vert_{t=1}+\|c_2(t)+\langle r(t):\phi(t)\rangle_t\|_{L^1_t[1,\infty)}\nonumber\\\leq & 2m^2 +\|c_2(t)\|_{L^1_t[1,\infty)} +\|\langle r(t):\phi(t)\rangle_t\|_{L^1_t[1,\infty)}<\infty.
\end{align}
 Hence, 
 \begin{align}
  \lim\limits_{T\to \infty}\int_1^{T} c_1(t)dt  \text{ exists and } \int_1^{\infty} c_1(t)dt<\infty.
 \end{align}
By using H\"{o}lder's inequality in $s$ variable, we have the following pointwise estimate: for $t_2\geq t_1\geq T\geq1$,
\begin{align}
|\psi_{\Omega,2}(t_2)&-\psi_{\Omega,2}(t_1)|=\int_{t_1}^{t_2}\partial_s[F_c]|F_1e^{isH_0}\psi(s)|ds\nonumber\\
&\leq\left(\int_{t_1}^{t_2}\partial_s[F_c]ds\right)^{1/2}\left(\int_{t_1}^{t_2}\partial_s[F_c]|F_1e^{isH_0}\psi(s)|^2ds\right)^{1/2}\\
&\leq\left(\int_{t_1}^{t_2}\partial_s[F_c]|F_1e^{isH_0}\psi(s)|^2ds\right)^{1/2}.
\end{align}
Taking $L^2_x$ norm and applying Fubini's theorem, we have that, for $t_2\geq t_1\geq T\geq 1$,
\begin{align}
	\|\psi_{\Omega,2}(t_2)-\psi_{\Omega,2}(t_1)\|_{L_x^2(\mathbb{R}^n)}\leq& \left(\int_{t_1}^{t_2}\|\sqrt{|\partial_s[F_c]|}F_1e^{isH_0}\psi(s)\|_{L_x^2(\mathbb{R}^n)}^2ds\right)^{1/2}\rightarrow 0,
\end{align}
as $T\rightarrow\infty.$ 
Hence $\{\psi_{\Omega,2}(t)\}_{t\geq 1}$ is Cauchy in $L^2_x(\R^n)$.
Therefore,
\begin{align}\label{limit1}
	\psi_{\Omega,2}(\infty):=\lim_{t\rightarrow\infty}\psi_{\Omega,2}(t)\ \mbox{exists in }\ L_x^2(\R^n).
\end{align}

For $\psi_{\Omega,3}(t)$,  we write $\psi_{\Omega,3}(t)$ as 
\begin{align}
\psi_{\Omega,3}(t)&=\int_{1}^t\partial_s[F_1]F_ce^{isH_0}\psi(s)ds-\int_{1}^t[\partial_s[F_1],F_c]e^{isH_0}\psi(s)ds:=\psi_{\Omega,31}(t)+\psi_{\Omega,32}(t).
\end{align}
For $\psi_{\Omega,32}(t)$, we have, by Corollary \ref{Cor1} and the condition that $b<\alpha\leq 1$,
\begin{equation}
    \|[\partial_s[F_1],F_c]e^{isH_0}\psi(s) \|\in L^1_s[1,\infty),
\end{equation}
and therefore conclude that
\begin{align}\label{limit2}
	\psi_{\Omega,32}(\infty):=\lim_{t\rightarrow\infty}\psi_{\Omega,32}(t)\text{ exists in $L^2_x(\mathbb{R}^n)$.}
\end{align}
For $\psi_{\Omega,31}(t)$, we use a similar argument as what we did for $\psi_{\Omega,2}(t)$.
 By taking 
$$
\begin{cases}
B(t)=F_cF_1F_c,\\
\phi(t)=e^{itH_0}\psi(t),
\end{cases}
$$ 
we have
\eq
\langle B(t): \phi(t)\rangle_t\leq \left(\sup\limits_{t\in \R}\|\psi(t)\|_{L^2_x(\mathbb{R}^n)}\right)^2=m^2<\infty.
\eeq
Since $$
\partial_t[F_1(t^b|p|>1)]=\frac{b}{t}\times \left(F'_1(t^b|p|>1)t^b|p|\right)\geq 0, 
$$
we compute
\eq
\partial_t[\langle B(t):\phi(t)\rangle_t]=\tilde{c}_{1}(t)+\tilde{c}_{2}(t)+\langle \tilde r(t):\phi(t)\rangle_t,
\eeq
where
\begin{align*}
\tilde{c}_{1}(t):=&(F_c \phi(t), \partial_t[F_1]F_c \phi(t))_{L^2_x(\R^n)}\geq 0,\quad \forall\, t\geq 1,
\end{align*}
\begin{align*}
\tilde{c}_2(t)&:=(-i)(F_ce^{itH_0}\psi(t), F_1F_ce^{itH_0} V(x,t)\psi(t))_{L^2_x(\R^n)}\\
&\quad+i(F_ce^{itH_0}V(x,t)\psi(t), F_1F_ce^{itH_0}\psi(t))_{L^2_x(\R^n)} 
\end{align*}
and  $\langle \tilde r(t):\phi(t)\rangle_t$ is defined in \eqref{r-tilde}.
Using  H\"{o}lder's inequality,  Proposition \ref{prop2.2} and Lemma \ref{lem2.5}, we get 
 \begin{align}
    |\tilde{c}_2(t)|\lesssim \frac{m}{t^{(\sigma+1)/2}}\|V(x,t)\psi(t)\|_{L^{\infty}_t L^2_{x,\sigma}(\mathbb{R}^{n+1})}\in L^1_t[1,\infty)
 \end{align}
 and 
 \begin{align}
     |\langle \tilde{r}(t):\phi(t)\rangle_t|\lesssim\frac{1}{t^{1+\alpha-b}}m^2\in L^1_t[1,\infty),
 \end{align}
 provided that $\sigma>1$ and $\alpha>b$.
 
 And again applying  \textbf{Relative Propagation Estimate}, we have, for all $T\geq 1,$
\begin{align}
\int_1^T \tilde{c}_{1}(t)dt\leq&\,  \langle B(t):\phi(t)\rangle_t\vert_{t=T}-\langle B(t):\phi(t)\rangle_t\vert_{t=1}+\|\tilde{c}_2(t)+\langle \tilde{r}(t):\phi(t)\rangle_t\|_{L^1_t[1,\infty)}\nonumber\\\leq &\, 2m^2 +\|\tilde{c}_2(t)\|_{L^1_t[1,\infty)}+\|\langle \tilde{r}(t):\phi(t)\rangle_t\|_{L^1_t[1,\infty)}   <\infty.
\end{align}
Hence, 
\eq
\int_1^\infty \tilde{c}_{1}(t)dt=\lim\limits_{T\to \infty}\int_1^T \tilde{c}_{1}(t)dt\quad \text{ exists in }\R.
\eeq

 By H\"{o}lder's inequality in $s$ variable, we have the following pointwise estimate:  for  $t_2\geq t_1\geq T>1,$
\begin{align}
|\psi_{\Omega,31}(t_2)&-\psi_{\Omega,31}(t_1)|=\int_{t_1}^{t_2}\partial_s[F_1]|F_ce^{isH_0}\psi(s)|ds\nonumber\\
&\leq\left(\int_{t_1}^{t_2}\partial_s[F_1]ds\right)^{1/2}\left(\int_{t_1}^{t_2}\partial_s[F_1]|F_ce^{isH_0}\psi(s)|^2ds\right)^{1/2}\nonumber\\
&\leq\left(\int_{t_1}^{t_2}\partial_s[F_1]|F_ce^{isH_0}\psi(s)|^2ds\right)^{1/2}.\label{eq: eq4}
\end{align}
We take $L^2_x(\mathbb{R}^n)$ norm for both sides of~\eqref{eq: eq4} and apply Fubini's theorem to obtain
\begin{align}
	&\|\psi_{\Omega,31}(t_2)-\psi_{\Omega,31}(t_1)\|_{L_x^2(\mathbb{R}^n)}\leq \left(\int_{t_1}^{t_2}\|\sqrt{\partial_s[F_1]}F_ce^{isH_0}\psi(s)\|_{L_x^2(\mathbb{R}^n)}^2ds\right)^{1/2}\rightarrow0,
\end{align}
as $T\rightarrow\infty.$ 
Hence $\{\psi_{\Omega,31}(t)\}_{t\geq 1}$ is Cauchy in $L^2_x(\R^n)$.
Therefore, we conclude that
\begin{align}\label{limit3}
	\psi_{\Omega,31}(\infty):=\lim_{t\rightarrow\infty}\psi_{\Omega,31}(t)\ \mbox{exists in }\ L_x^2(\R^n).
\end{align}

Combining  \eqref{limit1}, \eqref{limit2}, \eqref{limit3} and \eqref{limit4}, we obtain the existence of $\Omega_{\alpha,b}^*\psi_0$ in $L_x^2(\R^n).$

Next we prove \eqref{wlimit}. Take $g\in L^2_x(\R^n)$. We note that by Assumption~\ref{asp:sol},
\begin{align}
    &\left|(g,\left(1-F_c(\frac{|x|}{t^{\alpha}}\leq1)F_1(t^b|p|>1)\right)e^{itH_0}\psi(t))_{L^2_x(\mathbb{R}^n)}\right|\nonumber\\
    =&\left|(\left(1-F_c(\frac{|x|}{t^{\alpha}}\leq1)F_1(t^b|p|>1)\right)g,e^{itH_0}\psi(t))_{L^2_x(\mathbb{R}^n)}\right|\nonumber\\
    \leq & \|\left(1-F_c(\frac{|x|}{t^{\alpha}}\leq1)F_1(t^b|p|>1)\right)g \|_{L^2_x(\mathbb R^n)} m\to 0,
\end{align}
as $t\to \infty$. Thus the weak limit \eqref{wlimit} holds, which completes the proof.\end{proof}
\section{Localization of the weakly localized part}
In this section, we prove the localization properties of the non-free part, Theorem~\ref{Thm1} parts \emph{(ii)} and \emph{(iii)}. We need the notion of forward and backward waves. We include their definitions and properties in Section~\ref{sec: inout}.

\subsection{Forward/backward propagation waves}\label{sec: inout}
We start by discussing the concept of forward and backward propagation waves. These waves are analogous to the incoming and outgoing waves first introduced by Enss \cite{Enss1978}. 

Let $S^{n-1}$ denote the unit sphere in $\mathbb{R}^n$. We define a class of functions on $S^{n-1}$, $\{F^{\hat{h}}(\xi)\}_{\hat{h}\in I}$, as a smooth partition of unity with an index set 
 \eq
 I=\{ \hat{h}_1,\cdots, \hat{h}_N\}\subseteq S^{n-1}\label{indexI}
 \eeq
 for some $N\in \N^+$, satisfying that there exists $c>0$ such that for every $\hat{h}_i\in I$, 
 \eq
 F^{\hat{h}}(\xi)=\begin{cases}1 & \text{ when }|\xi-\hat{h}|<c\\ 0 & \text{ when }|\xi-\hat{h}|>2c\end{cases}, \quad \xi\in S^{n-1}.\label{ceq1}
 \eeq
Given $\hat{h}\in I$, we define $\tilde F^{\hat h}: S^{n-1}\to \R,$ as another smooth cut-off function satisfying 
 \eq
 \tilde{F}^{\hat{h}}(\xi)=\begin{cases}1 & \text{ when }|\xi-\hat{h}|<4c\\ 0 & \text{ when }|\xi-\hat{h}|>8c\end{cases},\quad \xi\in S^{n-1}.\label{ceq2}
 \eeq
For $h\in \mathbb{R}^n-\{0\}$, we define $\hat h:=h/|h|$ and $\hat h=0$ when $h=0$. We also assume that $c>0$, defined in \eqref{ceq1} and \eqref{ceq2}, is properly chosen  such that for all $x,q\in \R^n$ with $x\neq 0$ and $q\neq 0$, 
\eq
F^{\hat{h}}(\hat{x})\tilde{F}^{\hat{h}}(\hat{q})|x+q|\geq \frac{1}{10}(|x|+|q|)\label{Feq1}
\eeq
and
\eq
F^{\hat{h}}(\hat{x})(1-\tilde{F}^{\hat{h}}(\hat{q}))|x-q|\geq \frac{1}{10^6}(|x|+|q|).\label{Feq2}
\eeq
Now let us define the projection on forward/backward propagation set in terms of the phase-space $(r,v)\in \R^{n+n}$ with $|r|,|v|\neq 0$: 
\begin{definition}[Projection on the forward/backward propagation set] The projections onto the forward and backward propagation sets, in terms of $(r,v)\in \mathbb{R}^{n+n}$ with $|r|,|v|\neq 0$, are defined as follows:
\eq
P^+(r,v):=\sum\limits_{b=1}^N F^{\hat{h}_b}(\hat{r})\tilde{F}^{\hat{h}_b}(\hat{v}),\label{Prv+}
\eeq
and 
\eq
P^-(r,v):=1-P^+(r,v),\label{Prv-}
\eeq
respectively. 
\end{definition}

In our context, we take 
\begin{equation}\label{Ppm}
    P^\pm \equiv P^\pm (r,v)\vert_{r=x, v=4|p|^2p}.
\end{equation}
We need following estimates and their proofs can be found in Appendix.
\begin{lemma}\label{Lem: Pprop} For all $\epsilon \in (0,1/2)$ and $s,t,\sigma\geq 0$, the following estimates hold:
\begin{align}
&\| F_c(\frac{|x|}{(t+1)^{\frac{1}{4}+\epsilon}}\geq 1)P^\pm e^{\pm isH_0}F_1( (t+1)^{\frac{1}{4}}|p|\geq 1) \langle x\rangle^{-\sigma}\|_{L^2_x(\mathbb{R}^n)\to L^2_x(\mathbb{R}^n)}\nonumber\\
\lesssim_\epsilon &\frac{1}{\langle (t+1)^{\frac{1}{4}+\epsilon}+s/(t+1)^{3/4}\rangle^\sigma}\label{free: est: 1}
\end{align}
and
{\begin{align}
\| P^\pm e^{\pm isH_0} F_1((s+1)^{\frac{1}{4}-\epsilon}|p|\geq 1)\langle x\rangle^{-\sigma}\|_{L^2_x(\mathbb{R}^n)\to L^2_x(\mathbb{R}^n)}\lesssim_\epsilon \frac{1}{\langle s\rangle^{\sigma/4}}\label{free: est: 2}.
\end{align}}

\end{lemma}

\begin{lemma}\label{lem: Ppropfree} For all $\epsilon\in (0,1/2)$, $t,\sigma\geq 0$ and $s\in [0,t]$, the following estimate holds:
\begin{align}
&\| F_c(\frac{|x|}{(t+1)^{\frac{1}{4}+\epsilon}}\geq 1) e^{ -isH_0}F_1( (t+1)^{\frac{1}{4}}|p|< 1) \langle x\rangle^{-\sigma}\|_{L^2_x(\mathbb{R}^n)\to L^2_x(\mathbb{R}^n)}\nonumber\\
\lesssim_\epsilon &\frac{1}{ (t+1)^{\frac{1}{4}\sigma +\epsilon\sigma}}.\label{free: est: 3}
\end{align}

\end{lemma}

\begin{lemma}\label{lem: Ppmf} For all $f\in L^2_x(\mathbb{R}^n)$, we have
\begin{equation}
    \lim\limits_{s\to \infty}\| P^\pm e^{\pm i s H_0}f \|_{L^2_x(\mathbb{R}^n)}=0.
\end{equation}
    
\end{lemma}

\begin{lemma}\label{lem: not linear: local}
For all $f\in L^2_x(\mathbb{R}^n)$, $\alpha\in (0,1)$ and $s\geq 0$, we have
\begin{equation}
    \lim\limits_{s\to \infty}\|\chi(|x|\leq s^\alpha) P^\mp e^{\pm i s H_0}f \|_{L^2_x(\mathbb{R}^n)}=0.\label{goal: lem: eq1}
\end{equation} 
\end{lemma}

\subsection{Intertwining properties}
When the interaction~$V(x,t)=W(x)$ is real, time-independent and short-range, with $H=H_0+W(x)$, define 
\begin{equation}
\Omega_\pm:=s\text{-}\lim\limits_{t\to \infty} e^{itH}e^{-itH_0},\qquad \text{ on }L^2_x(\mathbb R^n)
\end{equation}
and
\begin{equation}
\Omega_\pm^*:=s\text{-}\lim\limits_{t\to \infty} e^{itH_0}e^{-itH}P_c,\qquad \text{ on }L^2_x(\mathbb R^n),
\end{equation}
where $P_c$ denotes the projection onto the continuous spectrum of $H$. These operators satisfy the intertwining relations
\begin{equation}
\Omega_\pm e^{-itH_0}=e^{-itH}\Omega_\pm \qquad \text{ on }L^2_x(\mathbb R^n)
\end{equation}
and
\begin{equation}
\Omega_\pm^* e^{-itH}=e^{-itH_0}\Omega_\pm^*\qquad \text{ on }L^2_x(\mathbb R^n).
\end{equation}

In this paper, introduce \emph{adapted wave operators} to obtain an adapted intertwining property for a broader class of bi-Laplacian Schr\"odinger equations. This framework will be used throughout the proofs in Sections~\ref{sec: 4.3} and~\ref{sec: 5}. 

More precisely, fix $\sigma>1$. For $t\in \mathbb R$, $b\in (0,(1/\sigma+1)/6)$ and $\alpha \in (0,(1+1/\sigma)/2)$, set
\begin{equation}
\Omega_{\alpha,b}(t)^*\psi(t):=s\text{-}\lim\limits_{s\to \infty} F_c(\frac{|x|}{s^\alpha}\leq 1)F_1(s^b|p|>1) e^{isH_0}\psi(t+s),\qquad \text{ in }L^2_x(\mathbb R^n).
\end{equation}
The existence of $\Omega_{\alpha,b}(t)^*\psi(t)$ in $L^2_x(\mathbb R^n)$ follows by the same argument as in Theorem~\ref{Thm1} part (i), provided that Assumptions~\ref{asp:sol} and~\ref{condition} hold. In what follows, we adopt the shorthand $F_c\equiv F_c(\frac{|x|}{s^\alpha}\leq 1)$ and $F_1\equiv F_1(s^b|p|>1)$. 

\begin{proposition}[Adapted intertwining identity]\label{prop: inter} Let $\sigma$ be as in Assumption~\ref{condition}. If Assumptions~\ref{asp:sol} and~\eqref{condition} hold, then for each $t\in \mathbb R$, $b\in (0,(1/\sigma+1)/6)$ and $\alpha \in (0,(1+1/\sigma)/2)$, we have 
\begin{equation}\label{eq: inter 1}
    \Omega_{\alpha,b}^*(t)\psi(t)=e^{-itH_0}\Omega_{\alpha,b}^*(0)\psi_0.
\end{equation}
Furthermore, define
\begin{equation}    \Omega^*(t)\psi(t):=w\text{-}\lim\limits_{s\to \infty} e^{isH_0}\psi(t+s),\qquad \text{ in }L^2_x(\mathbb R^n).
\end{equation}
Then $\Omega^*(t)\psi(t)$ exists in $L^2_x(\mathbb R^n)$, and 
\begin{equation}\label{eq: inter 2}
    \Omega^*(t)\psi(t)=\Omega_{\alpha,b}^*(t)\psi(t)
\end{equation}
and
\begin{equation}\label{eq: inter 3}
     \Omega^*(t)\psi(t)=e^{-itH_0}\Omega^*(0)\psi_0
\end{equation}
hold in the weak topology of $L^2_x(\mathbb R^n)$.
\end{proposition}
\begin{proof} We first prove the existence of $\Omega^*(t)\psi(t)$ in $L^2_x(\mathbb R^n)$. Fix $g\in L^2_x(\mathbb R^n)$. Since 
\begin{equation}
    s\text{-}\lim\limits_{s\to \infty} (1-F_1(s^b|p|>1)F_c(\frac{|x|}{s^\alpha}\leq 1))=0 \quad \text{ on }L^2_x(\mathbb R^n)
\end{equation}
for all $\alpha,b\in (0,1)$,  Assumption~\ref{asp:sol} implies,
\begin{equation}
    \begin{aligned}
      \limsup\limits_{s\to \infty}  \left|(g, e^{isH_0}\psi(s+t))_{L^2_x(\mathbb R^n)}-(F_1(s^b|p|>1)F_c(\frac{|x|}{s^\alpha}\leq 1)g, e^{isH_0}\psi(s+t))_{L^2_x(\mathbb R^n)}\right|=0,
    \end{aligned}
\end{equation}
which is equivalent to 
\begin{equation}
    \begin{aligned}
      \limsup\limits_{s\to \infty}  \left|(g, e^{isH_0}\psi(s+t))_{L^2_x(\mathbb R^n)}-(g, F_c(\frac{|x|}{s^\alpha}\leq 1)F_1(s^b|p|>1)e^{isH_0}\psi(s+t))_{L^2_x(\mathbb R^n)}\right|=0.
    \end{aligned}
\end{equation}
Together with Theorem~\ref{Thm1} part (i) yields 
    \begin{equation}
    \begin{aligned}
      \limsup\limits_{s\to \infty}  \left|(g, e^{isH_0}\psi(s+t))_{L^2_x(\mathbb R^n)}-(g, \Omega_{\alpha,b}^*(t)\psi(t))_{L^2_x(\mathbb R^n)}\right|=0.
    \end{aligned}
\end{equation}
Hence,  $\Omega^*(t)\psi(t)$ exists in $L^2_x(\mathbb R^n)$ in the weak topology,and consequently ~\eqref{eq: inter 2} and~\eqref{eq: inter 3} hold. 

We now prove~\eqref{eq: inter 1}. By~\eqref{eq: inter 2} and~\eqref{eq: inter 3}, 
\begin{equation}
    \begin{aligned}
        \Omega_{\alpha,b}^*(t)\psi(t)=\Omega^*(t)\psi(t)=e^{-itH_0}\Omega^*(0)\psi(0)=e^{-itH_0}\Omega_{\alpha,b}^*(0)\psi_0
        \end{aligned}
\end{equation}
in the weak topology. Additonally, by Theorem~\ref{Thm1} part (i),  both $\Omega_{\alpha,b}^*(t)\psi(t)$ and $\Omega_{\alpha,b}^*(0)\psi(0)$ exist in $L^2_x(\mathbb R^n)$ in the strong topology. Therefore  \eqref{eq: inter 1} holds  in the strong topology, which completes the proof.\end{proof}

\subsection{Proof of Theorem~\ref{Thm1} part \emph{(ii)}}\label{sec: 4.3}
We first outline the strategy for proving the weakly localization of the non-free component, using several auxiliary estimates established in Section~\ref{sec: Auxiliary} below.  Throughout, we adopt the convention $F_1\equiv F_1((t+1)^{\frac{1}{4}}|p|\geq 1)$ whenever the context is clear.\par 

To prove~\eqref{wlocalize2}, it suffices to show
\begin{equation}\label{eq: eq: 4.19}
    \| F_c( \frac{|x|}{(t+1)^{\frac{1}{4}+\epsilon}}\geq 1)\psi(t)-e^{-itH_0}\Omega_{\alpha,b}^*\psi_0\|_{L^2_x(\mathbb R^n)}\to 0
\end{equation}
as $t\to \infty$. For this purpose,  decompose
\begin{align}
F_c( \frac{|x|}{(t+1)^{\frac{1}{4}+\epsilon}}\geq 1)\psi(t)= & \psi_1(t)+\psi_2(t)+\psi_3(t),\label{decom: psi largex}
\end{align}
where $\psi_j(t), j=1,2,3,$ are defined by 
\eq
\psi_1(t):=F_c( \frac{|x|}{(t+1)^{\frac{1}{4}+\epsilon}}\geq 1)P^-F_1 \psi(t),
\eeq
\eq
\psi_2(t):=F_c( \frac{|x|}{(t+1)^{\frac{1}{4}+\epsilon}}\geq 1)(1- F_1)\psi(t)
\eeq
and
\eq
\psi_3(t):=F_c( \frac{|x|}{(t+1)^{\frac{1}{4}+\epsilon}}\geq 1)P^+F_1 \psi(t).
\eeq
We approximate $\psi_3(t)$ by $e^{-itH_0}\Omega^*_{\alpha,\beta}\psi_0$ and obtain ~\eqref{eq: eq: 4.19} from
\eq
\lim\limits_{t\to \infty}\| \psi_j(t)\|_{L^2_x(\mathbb{R}^n)}=0,\label{lim: psi2}\quad j=1,2,
\eeq
and
\begin{equation}\label{eq: eq: psi3}
    \lim\limits_{t\to \infty} \|\psi_3(t)-e^{-itH_0}\Omega_{\alpha,b}^*\psi_0\|_{L^2_x(\mathbb R^n)}=0,
\end{equation}
whose proofs are given in Section~\ref{sec: Auxiliary}.

\textbf{Idea behind~\eqref{lim: psi2} and~\eqref{eq: eq: psi3}.} By Duhamel's formula together with Lemmas~\ref{Lem: Pprop}, \ref{lem: Ppropfree}, and~\ref{lem: Ppmf}, the limit~\eqref{lim: psi2} follows by an argument similar to the proof of~\eqref{eq: eq: psi3}, presented below. To prove  ~\eqref{eq: eq: psi3}, by Proposition~\ref{prop: inter},  write

 \begin{equation}
     \begin{aligned}
    \psi_3(t)-e^{-itH_0} \Omega^*_{\alpha,b}\psi_0& =\left(F_c( \frac{|x|}{(t+1)^{\frac{1}{4}+\epsilon}}\geq1)P^+e^{-itH_0}F_1\Omega^*_{\alpha,b}\psi_0-e^{-itH_0} \Omega^*_{\alpha,b}\psi_0\right)\\
    &\quad+\left(F_c( \frac{|x|}{(t+1)^{\frac{1}{4}+\epsilon}}\geq 1)P^+F_1(1-\Omega^*_{\alpha,b}(t))\psi(t)\right) \\
    &=: \psi_{3,1}(t)+\psi_{3,2}(t).
     \end{aligned}
 \end{equation}
By Lemmas~\ref{lem: Ppmf} and~\ref{lem: not linear: local}, we have 
 \begin{equation}
   \begin{aligned}
    & \| \psi_{3,1}(t)\|_{L^2_x(\mathbb{R}^n)} \\
    \leq & \| F_c( \frac{|x|}{(t+1)^{\frac{1}{4}+\epsilon}}\geq 1)P^+(1-F_1) e^{-itH_0} \Omega^*_{\alpha,\beta}\psi(0) \|_{L^2_x(\mathbb{R}^n)} \\
    & +\| F_c( \frac{|x|}{(t+1)^{\frac{1}{4}+\epsilon}}< 1)P^+e^{-itH_0} \Omega^*_{\alpha,\beta}\psi(0) \|_{L^2_x(\mathbb{R}^n)}+\| P^-e^{-itH_0} \Omega^*_{\alpha,\beta}\psi(0) \|_{L^2_x(\mathbb{R}^n)} \\
    \leq & \|(1-F_1) \Omega^*_{\alpha,\beta}\psi(0) \|_{L^2_x(\mathbb{R}^n)}+\| F_c( \frac{|x|}{(t+1)^{\frac{1}{4}+\epsilon}}< 1)P^+e^{-itH_0} \Omega^*_{\alpha,\beta}\psi(0) \|_{L^2_x(\mathbb{R}^n)} \\
    & +\| P^-e^{-itH_0} \Omega^*_{\alpha,\beta}\psi(0) \|_{L^2_x(\mathbb{R}^n)}\to  0,\label{psi(31): local}
\end{aligned} 
\end{equation}
as $t\to \infty.$ Applying Duhamel’s formula and Lemma~\ref{Lem: Pprop} yields 
 \begin{align}\label{lim:psi3j}
     \lim\limits_{t\to \infty}\| \psi_{3,2}(t)\|_{L^2_x(\mathbb{R}^n}=0.
 \end{align}
Combining~\eqref{psi(31): local} and~\eqref{lim:psi3j} gives~\eqref{eq: eq: psi3}.

\subsection{Auxiliary estimates}\label{sec: Auxiliary}

\begin{proof}[Proof of \eqref{lim: psi2}.] By Duhamel's formula, $\psi_1(t)$ and $\psi_2(t)$ can be written as
\begin{align}
    \psi_1(t)=& F_c (\frac{|x|}{(t+1)^{\frac{1}{4}+\epsilon}}\geq 1) P^-F_1e^{-itH_0}\psi(0)\nonumber\\
    &+(-i)\int_0^t F_c (\frac{|x|}{(t+1)^{\frac{1}{4}+\epsilon}}\geq 1) P^-F_1e^{-i(t-s)H_0}V(x,s) \psi(s)ds\nonumber\\
    =: & \psi_{11}(t)+\psi_{12}(t)\label{psi21psi22}
\end{align}
and
\begin{align}
    \psi_2(t)=& F_c (\frac{|x|}{(t+1)^{\frac{1}{4}+\epsilon}}\geq 1) (1-F_1)e^{-itH_0}\psi(0)\nonumber\\
    &+(-i)\int_0^t F_c (\frac{|x|}{(t+1)^{\frac{1}{4}+\epsilon}}\geq 1) (1-F_1)e^{-i(t-s)H_0}V(x,s) \psi(s)ds\nonumber\\
    =: & \psi_{21}(t)+\psi_{22}(t),\label{psi21psi22'}
\end{align}
respectively. By Lemma~\ref{lem: Ppmf}, $\psi_{11}(t)$ satisfies 
\begin{align}
\| \psi_{11}(t)\|_{L^2_x(\mathbb{R}^n)}\leq& \| P^-e^{-itH_0}\psi(0)\|_{L^2_x(\mathbb{R}^n)}\to 0,\label{psi21(t)}
\end{align}
as $t\to \infty$. Using Lemma~\ref{Lem: Pprop} and Assumptions~\ref{asp:sol} and ~\ref{condition}, $\psi_{12}(t)$ obeys, for $\sigma\geq 4$, 
\begin{equation}
    \begin{aligned}
   & \| \psi_{12}(t)\|_{L^2_x(\mathbb{R}^n)}\\\leq & \int_0^t \| F_c (\frac{|x|}{(t+1)^{\frac{1}{4}+\epsilon}}\geq 1) P^-e^{-i(t-s)H_0} F_1 \langle x\rangle^{-\sigma}\|_{L^2_x(\mathbb{R}^n)\to L^2_x(\mathbb{R}^n)} \\
    &\times \| \langle x\rangle^\sigma V(x,s)\psi(s)\|_{L^2_x(\mathbb{R}^n)}ds \\
    \lesssim&_\epsilon  \int_0^t \frac{1}{\langle  (t+1)^{\frac{1}{4}+\epsilon}+s/(t+1)^{3/4}\rangle^{\sigma}} \sup\limits_{u\in \mathbb{R}} \| \langle x\rangle^\sigma V(x,u)\psi(u)\|_{L^2_x(\mathbb{R}^n)}ds  \\
    \lesssim& _{\epsilon,\delta} \frac{t}{\langle t+1\rangle^{1+4\epsilon}} \sup\limits_{u\in \mathbb{R}}\| \langle x\rangle^\sigma V(x,u)\psi(u)\|_{L^2_x(\mathbb{R}^n)} 
    \to  0,\label{psi22}
\end{aligned}
\end{equation}
as $t\to \infty.$ The bounds \eqref{psi21(t)} and~\eqref{psi22}, together with Eq.~\eqref{psi21psi22}, imply 
\eq
\lim\limits_{t\to \infty} \| \psi_1(t)\|_{L^2_x(\mathbb{R}^n)}=0.\label{lim psi20}
\eeq
We also note that
\eq
s\text{-}\lim\limits_{t\to \infty}1- F_1 =0,\quad \text{ on }L^2_x(\mathbb{R}^n),\label{F1}
\eeq
which yields
\eq
\lim\limits_{t\to \infty} \| \psi_{21}(t)\|_{L^2_x(\mathbb{R}^n)}=0.\label{lim psi3}
\eeq
Using Lemma~\ref{lem: Ppropfree} and Assumptions~\ref{asp:sol} and ~\ref{condition}, $\psi_{22}(t)$ obeys, for $\sigma\geq 4$, 
\begin{equation}
    \begin{aligned}
   & \| \psi_{22}(t)\|_{L^2_x(\mathbb{R}^n)}\\\leq & \int_0^t \| F_c (\frac{|x|}{(t+1)^{\frac{1}{4}+\epsilon}}\geq 1) e^{-i(t-s)H_0}(1- F_1) \langle x\rangle^{-\sigma}\|_{L^2_x(\mathbb{R}^n)\to L^2_x(\mathbb{R}^n)} \\
    &\times \| \langle x\rangle^\sigma V(x,s)\psi(s)\|_{L^2_x(\mathbb{R}^n)}ds \\
    \lesssim&_\epsilon  \int_0^t \frac{1}{(1+t)^{\frac{1}{4}\sigma+\epsilon \sigma}} \sup\limits_{u\in \mathbb{R}} \| \langle x\rangle^\sigma V(x,u)\psi(u)\|_{L^2_x(\mathbb{R}^n)}ds  \\
    \lesssim& _{\epsilon,\delta} \frac{t}{\langle t+1\rangle^{1+4\epsilon}} \sup\limits_{u\in \mathbb{R}}\| \langle x\rangle^\sigma V(x,u)\psi(u)\|_{L^2_x(\mathbb{R}^n)} 
    \to  0,\label{psi22'}
\end{aligned}
\end{equation}
as $t\to \infty.$ The bounds \eqref{lim psi3} and~\eqref{psi22'}, together with Eq.~\eqref{psi21psi22'}, imply 
\eq
\lim\limits_{t\to \infty} \| \psi_2(t)\|_{L^2_x(\mathbb{R}^n)}=0.\label{lim psi20'}
\eeq

Combining \eqref{lim psi20} and \eqref{lim psi20'} proves \eqref{lim: psi2}.
\end{proof}

Next, we prove~\eqref{lim:psi3j}. 
\begin{proof}[Proof of~\eqref{lim:psi3j}]
By Lemma~\ref{Lem: Pprop} and Assumption~\ref{asp:sol},  with $\sigma\geq 4$, as $t\to \infty$,
\begin{equation}
    \begin{aligned}
       &\| \psi_{3,2}(t)\|_{L^2_x(\mathbb{R}^n)} \\
     \leq & \| F_c(\frac{|x|}{(t+1)^{\frac{1}{4}+\epsilon}}\geq 1) P^+\int_t^\infty F_1 e^{i(s-t)H_0}V(x,s)\psi(s)ds\|_{L^2_x(\mathbb R^n)}\\
   \leq&  \int_t^\infty\|F_c(\frac{|x|}{(t+1)^{\frac{1}{4}+\epsilon}}\geq 1) P^+ F_1 e^{i(s-t)H_0}\langle x\rangle^{-\sigma}\|_{L^2_x(\mathbb{R}^n)\to L^2_x(\mathbb{R}^n)} \\
    &\times \|\langle x\rangle^\sigma V(x,s)\psi(s)\|_{L^2_x(\mathbb{R}^n)}ds 
    \\
     \lesssim&_\epsilon\int_t^\infty \frac{1}{\langle (t+1)^{\frac{1}{4}+\epsilon}+(s-t)/(t+1)^{3/4}\rangle^\sigma}\|\langle x\rangle^\sigma V(x,s)\psi(s)\|_{L^2_x(\mathbb{R}^n)}ds \\
    \lesssim& _{\epsilon,\delta} \frac{(t+1)^{3/4}}{\langle t+1\rangle^{\frac{1}{4}(\sigma-1)+\epsilon(\sigma -1)}}\|\langle x\rangle^\sigma V(x,s)\psi(s)\|_{L^2_x(\mathbb{R}^n)}   \to 0.\label{psi(32): local}
    \end{aligned}
\end{equation}
Together with \eqref{psi(31): local}, this proves~\eqref{lim:psi3j}.   
\end{proof}
\subsection{Proof of Theorem~\ref{Thm1} part \emph{(ii)}}
\begin{proof}[Proof of Theorem~\ref{Thm1} part (ii)] Combining \eqref{lim: psi2}, \eqref{eq: eq: psi3}, and the decomposition \eqref{decom: psi largex} yields \eqref{eq: eq: 4.19}, and hence \eqref{wlocalize2}.
\end{proof}

\section{Full localization of the non-free part in higher space dimension}\label{sec: 5}
In this section we establish part \emph{(iii)} of Theorem~\ref{Thm1} using Proposition~\ref{prop: weight} stated below. Throughout, we adopt the convention $F_{1,k}\equiv F_1\big((1+s+2^k)^{\frac{1}{4}}|p|>1\big)$ for $k\in\mathbb N$. Let
\begin{equation}\label{eq: loc}
\begin{aligned}
\psi_{loc}(t):=& P^+(1-\Omega_{\alpha,b}^*(t))\psi(t)+P^-(\psi(t)-e^{-itH_0}\psi_0).
\end{aligned}
\end{equation}

\begin{proposition}\label{prop: weight}If the assumptions stated in Theorem~\ref{Thm1} part (iii) are satisfied, then for all $\delta \in [0, \frac{n}{8}-1)$ with $n\geq 9$,
    \begin{equation}
      \sup\limits_{t\in \mathbb R}  \|\langle x\rangle^\delta\psi_{loc}(t)\|\lesssim_{\delta,n, \delta+1-n/8} 1.\label{est: Pbw}
    \end{equation}
\end{proposition}
The proof of Proposition~\ref{prop: weight} relies on the following lemmas, which are established in the Appendix.
\begin{lemma}\label{lem: gkpsi} For all $s\geq 0$ and $k\in \mathbb N$, we have 
\begin{equation}\label{eq: est 5.2}
    \| g_k(|x|)P^\pm e^{\pm isH_0}F_{1,k}\langle x\rangle^{-\sigma}\|_{L^2_x(\mathbb R^n)\to L^2_x(\mathbb{R}^n)}\lesssim_{\sigma,n} \frac{1}{\langle 2^k+\frac{s}{(2^k+s+1)^{\frac{3}{4}}}\rangle^\sigma}.
\end{equation}
    
\end{lemma}

\begin{lemma}\label{lem: gkpsi'} For all $n\in \mathbb N^+$, $k\in \mathbb N$ and $s\geq 0$, we have
\begin{equation}
    \| 1-F_{1,k}\|_{L^1_x(\mathbb R^n)\to L^2_x(\mathbb R^n)}\lesssim_n \frac{1}{\langle1+s+2^k\rangle^{\frac n2}}.
\end{equation}

\end{lemma}

\begin{proof}[Proof of Proposition~\ref{prop: weight}] We use the Littlewood-Paley decomposition to write: 
\begin{equation}\label{eq: LP decom}
    \psi_{loc}(t)=\sum\limits_{k=0}^\infty g_k(|x|)\psi_{loc}(t)
\end{equation}
with a radial partition of unity $\{g_k(|x|)\}_{k\in \mathbb N}$ satisfying 
\begin{equation}
    \text{supp}\,g_k\subseteq [2^{k-3},2^{k+3}],\qquad \forall \, k=1,2,\cdots
\end{equation}
and
\begin{equation}
    \text{supp}\,g_0\subseteq [0,8].
\end{equation}
For each block $g_k(|x|)\psi_{loc}(t)$, we further decompose
\begin{equation}\label{eq: def loc kpm}
    g_k(|x|)\psi_{loc}(t)=\psi_{loc,k,+}(t)+\psi_{loc,k,-}(t)
\end{equation}
with
\begin{equation}
    \psi_{loc,k,+}(t):=g_k(|x|)  P^+(1-\Omega_{\alpha,b}^*(t))\psi(t)
\end{equation}
and
\begin{equation}
    \psi_{loc,k,-}(t):=g_k(|x|)  P^-(\psi(t)-e^{-itH_0}\psi_0).
\end{equation}

Next, we estimate
\begin{equation}
    \begin{aligned}
       & \|\langle x\rangle^\delta \psi_{loc,k,\pm}(t)\|_{L^2_x(\mathbb R^n)}\\
       \leq &\int_0^\infty \|\langle x\rangle^\delta g_k(|x|) P^\pm e^{\pm isH_0}F_{1,k}V(t+s)\psi(t+s)\|_{L^2_x(\mathbb R^n)} ds\\
        &+\int_0^\infty \|\langle x\rangle^\delta g_k(|x|) P^\pm e^{\pm isH_0}(1-F_{1,k})V(t+s)\psi(t+s)\|_{L^2_x(\mathbb R^n)} ds.
    \end{aligned}
\end{equation}
By Lemmas~\ref{lem: gkpsi} and~\ref{lem: gkpsi'} together with Duhamel’s principle, we obtain
\begin{equation}
    \begin{aligned}
     \|\langle x\rangle^\delta&\psi_{loc,k,\pm}(t)\|_{L^2_x(\mathbb R^n)}\lesssim_{\sigma,n}  \int_0^\infty 2^{k\delta} \frac{1}{\langle 2^k+\frac{s}{(2^k+s+1)^{\frac{3}{4}}}\rangle^\sigma} \|\langle x\rangle^{\sigma}V(t+s)\psi(t+s)\|_{L^2_x(\mathbb R^n)}ds\\
        &+\int_0^\infty 2^{k\delta} \frac{1}{(2^k+s+1)^{\frac{n}{8}}}  \|\langle x\rangle^{\sigma}V(t+s)\psi(t+s)\|_{L^2_x(\mathbb R^n)}ds.
        \end{aligned}
        \end{equation}
Combining this with the elementary bound
\begin{equation}
\begin{aligned}
    \frac{1}{\langle 2^k+\frac{s}{(2^k+s+1)^{\frac{3}{4}}}\rangle^\sigma}\lesssim & \frac{\chi(s\leq 2^k)}{2^{k\sigma}}+\frac{\chi(s>2^k)}{(s+1)^{\frac{\sigma}{4}}}\\
    \lesssim &\frac{\chi(s\leq 2^k)}{2^{kn/2}}+\frac{\chi(s>2^k)}{(s+1)^{\frac{n}{8}}}\\
    \lesssim& \frac{1}{( 2^k+s+1)^{\frac{n}{8}}}
    \end{aligned},\qquad \forall\, s\geq 0,\, k\in \mathbb N\text{ and }\sigma>\frac{n}{2},
\end{equation}
we conclude that
\begin{equation}
        \begin{aligned}
        \|\langle x\rangle^\delta&\psi_{loc,k,\pm}(t)\|_{L^2_x(\mathbb R^n)}\lesssim_{\sigma,n}  \int_0^\infty \frac{2^{k\delta}}{(2^k+s+1)^{\frac{n}{8}}}  \|\langle x\rangle^{\sigma}V(t+s)\psi(t+s)\|_{L^2_x(\mathbb R^n)}ds.
    \end{aligned}
\end{equation}
Summing in $k$ and using \eqref{eq: LP decom}–\eqref{eq: def loc kpm}, we obtain
\begin{equation}
    \begin{aligned}
        \|\langle x\rangle^\delta&\psi_{loc }(t)\|_{L^2_x(\mathbb R^n)} \lesssim_{\sigma,n}\sum_{k=0}^\infty  \int_0^\infty \frac{2^{k\delta}}{(2^k+s+1)^{\frac{n}{8}}} ds \|\langle x\rangle^{\sigma}V(t )\psi(t )\|_{L^\infty_tL^2_x(\mathbb R\times\mathbb R^n)}\\
        &\lesssim_{\sigma,n}\sum_{k=0}^{\infty}2^{k(\delta+1-\frac{n}{8} )} \|\langle x\rangle^{\sigma}V(t )\psi(t )\|_{L^\infty_tL^2_x(\mathbb R\times\mathbb R^n)}\\
        &\lesssim_{\sigma,n, \delta+1-n/8} \|\langle x\rangle^{\sigma}V(t )\psi(t )\|_{L^\infty_tL^2_x(\mathbb R\times\mathbb R^n)}.
    \end{aligned}
\end{equation}
where we used $\delta<\dfrac{n}{8}-1$ and $n\geq 9$ to ensure summability. Hence \eqref{est: Pbw} holds. This completes the proof.
\end{proof}
Next, we prove Theorem~\ref{Thm1} part \emph{(iii)}.
\begin{proof}[Proof of Theorem~\ref{Thm1} part (iii)] By~\eqref{wlocalize2} and Proposition~\ref{prop: weight}, it suffices to prove that 
\begin{equation}\label{eq: thm1.5iii goal}
    \|\psi(t)-e^{-itH_0}\Omega_{\alpha,b}^*\psi_0-\psi_{loc}(t) \|_{L^2_x(\mathbb R^n)}\to 0\qquad \text{ as }t\to \infty
\end{equation}
holds. For this, by~\eqref{eq: inter 1} and~\eqref{eq: inter 2}, we write 
\begin{equation}
    \psi(t)=P^+\psi(t)+P^-\psi(t)=P^+e^{-itH_0}\Omega_{\alpha,b}^*\psi_0+P^-e^{-itH_0}\psi_0+\psi_{loc}(t),
\end{equation}
that is,
\begin{equation}
    \psi(t)=e^{-itH_0}\Omega_{\alpha,b}^*\psi_0-P^-e^{-itH_0}\Omega_{\alpha,b}^*\psi_0+P^-e^{-itH_0}\psi_0+\psi_{loc}(t).
\end{equation}
This together with Lemma~\ref{lem: Ppmf} yields~\eqref{eq: thm1.5iii goal} which completes the proof.
\end{proof}

\appendix 
\section{Auxiliary estimates}
\begin{proof}[Proof of Lemma~\ref{Lem: Pprop}] In this proof, we use the shorthand
\begin{equation}
    F_c=F_c(\frac{|x|}{(t+1)^{\frac{1}{4}+\epsilon}}\geq 1) 
\end{equation}
and
\begin{equation}
    F_1=F_1((t+1)^{\frac{1}{4}}|p|\geq 1)\text{ or } F_1=F_1((t+1)^{\frac{1}{4}}|q|\geq 1).
\end{equation}
We begin with the proof of estimate~\eqref{free: est: 1}. Set
\eq
A^\pm(t,s):=  P^\pm e^{\pm isH_0}F_1\langle x\rangle^{-\sigma}.
\eeq
Decompose $A^\pm (t,s)$ into two two parts:
\begin{align}
    A^\pm(t,s)=&A^\pm_1(t,s)+A^\pm_2(t,s),
\end{align}
where 
\begin{align}
A^\pm_1(t,s):=& P^\pm e^{\pm isH_0}F_1 \langle x\rangle^{-\sigma}\chi(|x|> \frac{1}{10^{10}}((t+1)^{\frac{1}{4}+\epsilon}+s/(t+1)^{3/4}))
\end{align}
and
\begin{align}
A^\pm_2(t,s):=& P^\pm e^{\pm isH_0}F_1 \langle x\rangle^{-\sigma}\chi(|x|\leq \frac{1}{10^{10}}((t+1)^{\frac{1}{4}+\epsilon}+s/(t+1)^{3/4})).
\end{align}
For $\|A^\pm_1(t,s)\|_{L^2_x(\mathbb{R}^n)\to L^2_x(\mathbb{R}^n)}$, we have
\begin{align}
\|A^\pm_1(t,s)\|_{L^2_x(\mathbb{R}^n)\to L^2_x(\mathbb{R}^n)}\leq &\| \langle x\rangle^{-\sigma} \chi(|x|>\frac{1}{10^{10}}((t+1)^{\frac{1}{4}+\epsilon}+s/(t+1)^{3/4}))\|_{L^2_x(\mathbb{R}^n)\to L^2_x(\mathbb{R}^n)}\nonumber\\
\lesssim & \frac{1}{\langle (t+1)^{\frac{1}{4}+\epsilon}+s/(t+1)^{3/4}\rangle^\sigma}.\label{Apm1}
\end{align}
For $A^\pm_2(t,s)$, choose $f\in L^2_x(\mathbb{R}^n)$ and by taking the Fourier transform, we can write
\begin{align}
    A^\pm_2(t,s)f= & \frac{1}{(2\pi)^n} \int  P^\pm(x,4|q|^2q) e^{i\phi(q)} F_1   \nonumber\\
    & \times \langle y\rangle^{-\sigma}\chi(|y|\leq \frac{1}{10^{10}}((t+1)^{\frac{1}{4}+\epsilon}+s/(t+1)^{3/4}))f(y)d^nyd^nq\label{Apm2ts},
\end{align}
where we have used that $P^\pm=P^\pm(x,4|q|^2q)$ and   note that the phase   is given by $\phi(q)=(x-y)\cdot q \pm s |q|^4$.
When 
\eq
|y|\leq \frac{1}{10^{10}}((t+1)^{\frac{1}{4}+\epsilon}+s/(t+1)^{3/4}),
\eeq
\eq
|x|\geq \frac{(t+1)^{\frac{1}{4}+\epsilon}}{2}
\eeq
and
\eq
|q|\geq \frac{1}{4(t+1)^{\frac{1}{4}}},
\eeq
using Eqs.~\eqref{Feq1}-\eqref{Ppm}, we obtain
\begin{align}
|\nabla_q[\phi(q)]|=&|x-y\pm 4s|q|^2q|\geq  \frac{1}{2}|x\pm 4s|q|^2q|-|y|\geq  \frac{1}{10^6}(|x|+4s|q|^3)-|y|\nonumber\\
\geq & \frac{1}{10^7}(|x|+4s|q|^3 ).\label{Ap: B:est1}
\end{align}
Choose an orthonormal basis $\{e_1,\cdots,e_n\}$ with $e_1:=\frac{x-y\pm 4s|q|^2q}{|x-y\pm 4s|q|^2q|}$, and denote $x_1:=x\cdot e_1, y_1:=x\cdot e_1$ and $q_1:=q\cdot e_1$. 
Then
\begin{align}
    \| F_cA^\pm_2(t,s)f\|_{L^2_x(\mathbb{R}^n)}\leq & \| \langle x\rangle^{-n}\|_{L^2_x(\mathbb{R}^n)}\| \langle x\rangle^n F_c A^\pm_2(t,s)f\|_{ L^\infty_x(\mathbb{R}^n)}\lesssim \| \langle x\rangle^n A^\pm_2(t,s)f\|_{ L^\infty_x(\mathbb{R}^n)} .\label{2toinfty}
\end{align}
By  ~\eqref{Ap: B:est1},  
\begin{align}
|\partial_{q_1}[\frac{1}{(x_1-y_1\pm 4s|q|^2 q_1)}]|=& |\frac{\mp 4s(|q|^2+2q_1^2)}{ (x_1-y_1\pm 4s|q|^2q_1)^2}|\nonumber\\
\lesssim  & \frac{4s(|q|^2+2q_1^2)}{(|x|+s|q|^3)^2}\lesssim  \frac{1}{|q|(|x|+s|q|^3)}.
\end{align}
This, together with estimates
\begin{align}
    | \partial_{q_1}[F_1((t+1)^{\frac{1}{4}}|q|\geq 1)]|= & \frac{(t+1)^{\frac{1}{4}}|q_1|}{|q|} |F_1'((t+1)^{\frac{1}{4}}|q|\geq 1)|\nonumber\\
    \lesssim  &\frac{1}{|q|}, 
\end{align}
\begin{align}
    | \partial^j_{q_1}[F_1((t+1)^{\frac{1}{4}}|q|\geq 1)]|\lesssim_j  &\frac{1}{|q|^j} 
\end{align}
and (recall that $P^\pm(r,v)$ is defined via $\hat{r}$ and $\hat{v}$. See Eqs.~\eqref{Prv+} and~\eqref{Prv-}.)
\begin{align}
    |\partial_{q_1}^j[P^\pm(x,4s|q|^2 q)]|\lesssim_j \frac{1}{|q|^j},\quad j=1,2,\cdots.
\end{align}
Hence, integrating by parts $N$ times in $q_1$,
\begin{align}
|\langle x\rangle^n A_2^\pm(t,s)f|\lesssim & \int \frac{\langle x\rangle^n \chi(|q|\geq \frac{1}{4(t+1)^{\frac{1}{4}}})}{ |q|^N\langle |x|+s|q|^3\rangle^N} \nonumber\\
&\times\langle y\rangle^{-\sigma}\chi(|y|\leq \frac{1}{10^{10}}((t+1)^{\frac{1}{4}+\epsilon}+s/(t+1)^{3/4}))|f(y)|d^nqd^ny .\label{Apm2}
\end{align}
Integrating in $q$ in \eqref{Apm2} and using,  for    $|q|\geq 1/(4(t+1)^{\frac{1}{4}})$,
\eq
\frac{1}{\langle |x|+s|q|^3\rangle}\lesssim \frac{1}{\langle |x|+s/(64(t+1)^{3/4})\rangle}
\eeq
and
\eq
\frac{\langle x\rangle^n}{\langle |x|+s/(64(t+1)^{3/4})\rangle^N}\lesssim  \frac{1}{\langle |x|+s/(64(t+1)^{3/4})\rangle^{N-n}},
\eeq
we get
\begin{equation}
    \begin{aligned}
    &\|\langle x\rangle^n F_cA_2^\pm(t,s)f\|_{L^\infty_x(\mathbb{R}^n)}\\\lesssim &\int \frac{\langle x\rangle^nF_c \chi(|q|\geq \frac{1}{4(t+1)^{\frac{1}{4}}})}{ |q|^N\langle |x|+s/((t+1)^{3/4})\rangle^N}  \\
&\times\langle y\rangle^{-\sigma}\chi(|y|\leq \frac{1}{10^{10}}((t+1)^{1/2+\epsilon}+s/(t+1)^{3/4}))|f(y)|d^nqd^ny \\
\lesssim & \int \frac{F_c 
 (t+1)^{\frac{N-n}{4}}}{ \langle |x|+s/(64(t+1)^{3/4})\rangle^{N-n}}  \\
&\times\langle y\rangle^{-\sigma}\chi(|y|\leq \frac{1}{10^{10}}((t+1)^{\frac{1}{4}+\epsilon}+s/(t+1)^{3/4}))|f(y)|d^ny.\label{ApB: est2}
\end{aligned}
\end{equation}
By Cauchy–Schwarz in ~\eqref{ApB: est2}, we obtain 
\begin{equation}
    \begin{aligned}
    & \|\langle x\rangle^n F_cA_2^\pm(t,s)f\|_{L^\infty_x(\mathbb{R}^n)} \\
    \lesssim & \frac{F_c 
 (t+1)^{\frac{N-n}{4}}}{ \langle |x|+s/(64(t+1)^{3/4})\rangle^{N-n}} \|f(y)\|_{L^2_y(\mathbb{R}^n)} \\
 &\qquad \times\| \langle y\rangle^{-\sigma}\chi(|y|\leq \frac{1}{10^{10}}((t+1)^{\frac{1}{4}+\epsilon}+s/(t+1)^{3/4}))\|_{L^2_y(\mathbb{R}^n)} \\
 \lesssim & \frac{ 
 (t+1)^{\frac{N-n}{4}}}{ \langle (t+1)^{\frac{1}{4}+\epsilon}+s/(64(t+1)^{3/4})\rangle^{N-3n/2}} \|f(y)\|_{L^2_y(\mathbb{R}^n)} \\
 \lesssim& \frac{1}{ \langle (t+1)^{1+4 \epsilon}+s/(64(t+1)^{3/4})\rangle^{\frac{4\epsilon N}{1+4\epsilon}- \frac{n}{2}-\frac{4\epsilon }{1+4\epsilon }n}} \|f(y)\|_{L^2_y(\mathbb{R}^n)} \\
 \lesssim & \frac{ 
 1}{ \langle (t+1)^{\frac{1}{4}+\epsilon}+s/((t+1)^{3/4})\rangle^{\sigma}} \|f(y)\|_{L^2_y(\mathbb{R}^n)}\label{estimate: A2pm}
\end{aligned}
\end{equation}
with $N = [\frac{1+4\epsilon}{4\epsilon}\sigma+n+\frac{1+4\epsilon}{8\epsilon}n]+1$. Combining \eqref{Apm1}, \eqref{2toinfty}, and \eqref{estimate: A2pm} gives \eqref{free: est: 1}. The proof of \eqref{free: est: 2} is analogous.
\end{proof}

\begin{proof}[Proof of Lemma~\ref{lem: Ppropfree}]Let $x$ and $y$ denote the position variables on the left- and right-hand sides, respectively. The group velocity is $\nabla_p[H_0]=4|p|^2p$. Let $q$ be the Fourier variable. When
$$
|x|\geq \frac{(t+1)^{\frac{1}{4}+\epsilon}}{2},\qquad s|q|\leq \frac{2t}{(t+1)^{\frac{1}{4}}},\qquad |y|\leq (t+1)^{\frac{1}{4}+\epsilon}/4,
$$
we have
$$
|x-y-4s|q|^2q|\geq |x|-|y|-4s|q|^3\geq (t+1)^{\frac14+\epsilon}/20.
$$
Therefore, by an argument analogous to Lemma~\ref{Lem: Pprop}, we obtain estimate~\eqref{free: est: 3}. This completes the proof.
\end{proof}

\begin{proof}[Proof of Lemma~\ref{lem: Ppmf}] We fix $s\geq 0$. For all $M\geq 1$ and $\epsilon\in (0,1/2)$,  $P^\pm e^{\pm isH_0}f$ satisfies 
\begin{equation}
    \begin{aligned}
    &\| P^\pm e^{\pm isH_0}f\|_{L^2_x(\mathbb{R}^n)} \\\leq & \| P^\pm e^{\pm isH_0}F_1 \chi(|x|\leq M) f \|_{L^2_x(\mathbb{R}^n)} +\| P^\pm e^{\pm isH_0}F_1 \chi(|x|> M) f \|_{L^2_x(\mathbb{R}^n)} \\
    &+\| P^\pm e^{\pm isH_0}(1-F_1 )f \|_{L^2_x(\mathbb{R}^n)} \\
    \leq & \|  P^\pm e^{\pm isH_0}F_1\langle x\rangle^{-1}\|_{L^2_x(\mathbb{R}^n)\to L^2_x(\mathbb{R}^n)}  \|\langle x\rangle\chi(|x|\leq M) f \|_{L^2_x(\mathbb{R}^n)}+\| \chi(|x|>M)f\|_{L^2_x(\mathbb{R}^n)} \\
    &+\| (1-F_1  )f \|_{L^2_x(\mathbb{R}^n)}.
\end{aligned}
\end{equation}
Choosing  $M=(1+s)^{1/100}$ and applying Lemma~\ref{Lem: Pprop}, we get
\begin{equation}
    \begin{aligned}
    &\| P^\pm e^{\pm isH_0}f\|_{L^2_x(\mathbb{R}^n)}\\
    \lesssim&_\epsilon  \frac{1}{\langle s\rangle^{1/2}} (1+s)^{1/100} \| f\|_{L^2_x(\mathbb{R}^n)} +\| \chi(|x|>(s+1)^{1/100})f\|_{L^2_x(\mathbb{R}^n)} \\
    &+\| (1-F_1  )f \|_{L^2_x(\mathbb{R}^n)} \\
    \to & 0\label{psi21}
\end{aligned}
\end{equation}
as $s\to \infty.$  \end{proof}

\begin{proof}[Proof of Lemma~\ref{lem: not linear: local}] Let $x$ and $y$ denote the position variables on the left-hand side and the right-hand side, respectively. The velocity is   $\nabla_p[H_0]=4|p|^2p$. We decompose $\chi(|x|\leq s^\alpha)P^\mp e^{\pm i sH_0}f$ into three pieces:
\begin{equation}
    \begin{aligned}
&\chi(|x|\leq s^\alpha)P^\mp e^{\pm i sH_0}f\\
=&\,\chi(|x|\leq s^\alpha)P^\mp e^{\pm i sH_0}F_2(s^{(1-\alpha)/300}|p|>1 )\chi(|x|\leq s^\alpha)f \\
+&\,\chi(|x|\leq s^\alpha)P^\mp e^{\pm i sH_0}F_2(s^{(1-\alpha)/300}|p|\leq1 )\chi(|x|\leq s^\alpha)f \\
+&\,\chi(|x|\leq s^\alpha)P^\mp e^{\pm i sH_0}\chi(|x|> s^\alpha)f \\
 =:&\, f_1(s)+f_2(s)+f_3(s).
\end{aligned}
\end{equation}
We have 
\begin{align}
    \limsup\limits_{s\to \infty}\|f_2(s)\|_{L^2_x(\mathbb{R}^n)}\leq&\limsup\limits_{s\to \infty} \| F_2(s^{(1-\alpha)/300}|p|\leq 1) \chi(|x|\leq s^\alpha)f\|_{L^2_x(\mathbb{R}^n)}\nonumber\\
    \leq& \limsup\limits_{s\to \infty}\|F_2(s^{(1-\alpha)/300}|p|\leq1 )f\|_{L^2_x(\mathbb{R}^n)}+\limsup\limits_{s\to \infty}\| \chi(|x|> s^\alpha)f\|_{L^2_x(\mathbb{R}^n)}\nonumber\\
    =& 0\label{free: f2}
\end{align}
and
\begin{align}
    \limsup\limits_{s\to \infty}\|f_3(s)\|_{L^2_x(\mathbb{R}^n)}\leq& \limsup\limits_{s\to \infty}\| \chi(|x|> s^\alpha)f\|_{L^2_x(\mathbb{R}^n)}\nonumber\\
    = & 0.\label{free: f3}
\end{align}
To estimate $f_1(s)$, we follow the same strategy as in Lemma~\ref{Lem: Pprop}. Taking the Fourier transform, $f_1(s)$ can be written as
\begin{align}
    f_1(s)= & \frac{1}{(2\pi)^n} \int  \chi(|x|\leq s^\alpha)e^{i\phi(q)} P^\pm(x,4|q|^2q) F_2(s^{(1-\alpha)/300}|q|>1)  \nonumber\\
    & \times\chi(|y|\leq s^\alpha)f(y)d^nyd^nq\label{fts},
\end{align}
where we  used  $P^\pm=P^\pm(x,4|q|^2q)$ (see Eq.~\eqref{Ppm}) and  note that the phase   is  
$\phi(q)=(x-y)\cdot q \pm |q|^4.$ 

When $|y|\leq s^\alpha$, $|x|\leq s^\alpha$ and $|q|\geq \frac{1}{2}s^{(\alpha-1)/300}$, by Eqs.~\eqref{Feq1}-~\eqref{Ppm} and estimate
\begin{align}
8s|q|^3\geq s^{(99+\alpha)/100}\geq s^\alpha,\quad s\geq 1,
\end{align}
we have
\begin{align}
|\nabla_q[\phi(q)]|=&|x-y\pm 4s |q|^2q|\geq  \frac{1}{2}|x\pm 4s |q|^2q|-|y|\geq  \frac{1}{10^6}(|x|+4s |q|^3)-|y|\nonumber\\
\geq & \frac{1}{10^7}(|x|+4s |q|^3 ).\label{Ap: B:est10}
\end{align}
Choose an orthonormal basis  $\{e_1,\cdots,e_n\}$ with $e_1:=\frac{x-y\pm 4s|q|^2q}{|x-y\pm 4s|q|^2q|}$, and denote $x_1:=x\cdot e_1, y_1:=x\cdot e_1$ and $q_1:=q\cdot e_1$. 
Then
\begin{align}
    \| f_1(s)\|_{L^2_x(\mathbb{R}^n)}\leq & \| \langle x\rangle^{-n}\|_{L^2_x(\mathbb{R}^n)}\| \langle x\rangle^n f_1(s)\|_{ L^\infty_x(\mathbb{R}^n)}\lesssim \| \langle x\rangle^n f_1(s)\|_{ L^\infty_x(\mathbb{R}^n)} .\label{B2toinfty}
\end{align}
By  ~\eqref{Ap: B:est1}, 
\begin{align}
|\partial_{q_1}[\frac{1}{(x_1-y_1\pm 4s|q|^2 q_1)}]|=& |\frac{\mp 4s(|q|^2+2q_1^2)}{ (x_1-y_1\pm 4s|q|^2q_1)^2}|\nonumber\\
\lesssim  & \frac{4s(|q|^2+2q_1^2)}{(|x|+s|q|^ 3)^2}\lesssim  \frac{1}{|q|(|x|+s|q|^3)}.
\end{align}
This, together with estimates
\begin{align}
    | \partial_{q_1}[F_1(s^{(1-\alpha)/300}|q|\geq 1)]|= & \frac{s^{(1-\alpha)/300}|q_1|}{|q|} |F_1'(s^{(1-\alpha)/100}|q|\geq 1)|\lesssim  \frac{1}{|q|}, 
\end{align}
\begin{align}
    | \partial^j_{q_1}[F_1(s^{(1-\alpha)/300}|q|\geq 1)]|= & \frac{s^{(1-\alpha)/300}|q_1|}{|q|} |F_1'(s^{(1-\alpha)/300}|q|\geq 1)|\lesssim_j  \frac{1}{|q|^j}
\end{align}
and (recall that $P^\pm(r,v)$ is defined in terms of $\hat{r}$ and $\hat{v}$. See Eqs.~\eqref{Prv+} and~\eqref{Prv-}.)
\begin{align}
    |\partial_{q_1}^j[P^\pm(x,4s|q|^2 q)]|\lesssim_j \frac{1}{|q|^j},\quad j=1,2,\cdots.
\end{align}
 Hence, integrating by parts $N$ times in $q_1$,
\begin{align}
|\langle x\rangle^n f_1(s)|\lesssim & \int \frac{\langle x\rangle^n \chi(|q|\geq \frac{1}{2}s^{(\alpha-1)/300})}{ |q|^N\langle |x|+s|q|^ 3\rangle^N} \chi(|y|\leq s^\alpha)|f(y)|d^nqd^ny .\label{BApm2}
\end{align}
Integrating in $q$ in  ~\eqref{BApm2} and using, for $|q|\geq \frac{1}{2}s^{(\alpha-1)/300}$,
\eq
\frac{1}{\langle |x|+s|q|^3 \rangle}\lesssim \frac{1}{\langle |x|+s^{(99+\alpha)/100}\rangle}
\eeq
and
\eq
\frac{\langle x\rangle^n}{\langle |x|+s^{(99+\alpha)/100}\rangle^{N}}\lesssim  \frac{1}{\langle s^{(99+\alpha)/100}\rangle^{N-n}},
\eeq
we obtain 
\begin{align}
    \|\langle x\rangle^n f_1(s)\|_{L^\infty_x(\mathbb{R}^n)}\lesssim &\int \frac{\langle x\rangle^n \chi(|q|\geq \frac{1}{2}s^{(\alpha-1)/300})}{ |q|^N\langle |x|+s|q|^3\rangle^N} \chi(|y|\leq s^\alpha)|f(y)|d^nqd^ny\nonumber\\
\lesssim & \int \frac{ 
 s^{ (\alpha-1)(N-n)/300}}{ \langle s^{(99+\alpha)/100}\rangle^{N-n}} \chi(|y|\leq s^\alpha)|f(y)|d^ny\nonumber\\
 \lesssim & \frac{1}{\langle s\rangle^{N-n}}\int  \chi(|y|\leq s^\alpha)|f(y)|d^ny.\label{ApB: est20}
\end{align}
By Cauchy–Schwarz i in  ~\eqref{ApB: est2}, we arrive at, as $s\to \infty$, 
\begin{align}
     \|\langle x\rangle^n f_1(s)\|_{L^\infty_x(\mathbb{R}^n)}\lesssim & \frac{1}{\langle s\rangle^{N-n}} \| f(y)\|_{L^2_y(\mathbb{R}^n)}\| \chi(|y|\leq s^\alpha)\|_{L^2_y(\mathbb{R}^n)}\nonumber\\
     \lesssim &\frac{1}{\langle s\rangle^{N-n-\alpha n/2}} \| f(y)\|_{L^2_y(\mathbb{R}^n)} \nonumber\\
     \to & 0
\end{align}
with $N=n+\alpha n/2+1$. Together with \eqref{free: f2} and \eqref{free: f3}, this implies \eqref{goal: lem: eq1}. 
\end{proof}

\begin{proof}[Proof of Lemma~\ref{lem: gkpsi}]
    \begin{equation}
        \begin{aligned}
            g_k(|x|)P^\pm e^{\pm isH_0}F_{1,k}\langle x\rangle^{-\sigma}=A_{1,k}(t,s)+A_{2,k}(t,s)
        \end{aligned}
    \end{equation}
where 
\begin{align}
    A_{1,k}(t,s):=g_k(|x|)P^\pm e^{\pm isH_0}F_{1,k}\langle x\rangle^{-\sigma}\chi(|x|>\frac{1}{10^{10}}(2^k+\frac{s}{(1+s+2^k)^{\frac{3}{4}}}))
\end{align}
and 
\begin{align}
    A_{2,k}(t,s):=g_k(|x|)P^\pm e^{\pm isH_0}F_{1,k}\langle x\rangle^{-\sigma}\chi(|x|\leq \frac{1}{10^{10}}(2^k+\frac{s}{(1+s+2^k)^{\frac{3}{4}}}))
\end{align}
For $\| A_{1,k}(t,s)\|_{L^2_x(\mathbb R^n)\to L^2_x(\mathbb R^n)}$, we simply use the weight to obtain 
\begin{equation}
    \begin{aligned}
        \| A_{1,k}(t,s)\|_{L^2_x(\mathbb R^n)\to L^2_x(\mathbb R^n)}&\leq \|\langle x\rangle^{-\sigma}\chi(|x|>\frac{1}{10^{10}}(2^k+\frac{s}{(1+s+2^k)^{\frac{3}{4}}}))\|_{L^2_x(\mathbb R^n)\to L^2_x(\mathbb R^n)}\\
        &\lesssim \frac{1}{\langle 2^k+\frac{s}{(2^k+s+1)^{\frac{3}{4}}}\rangle^\sigma}.
    \end{aligned}
\end{equation}
For$\| A_{2,k}(t,s)\|_{L^2_x(\mathbb R^n)\to L^2_x(\mathbb R^n)}$, we use a nonstationary phase argument in Fourier variables. The group velocity associated with $H_0$ is $\nabla_p[H_0]=4|p|^2p$. Let $q$ denote the variable in the Fourier space. When $|x|\geq  2^{k-3}$, $|q|>\frac{1}{(2^k+s+1)^{\frac{1}{4}}}$ and $|y|\leq  \frac{1}{10^{10}}(2^k+\frac{s}{(1+s+2^k)^{\frac{3}{4}}})$, 
\eq
|x-y-4s|q|^2q|\geq \frac{1}{10^{6}}(|x|+4s|q|^3-|y|)\geq \frac{1}{10^{7} }(2^k+4s|q|^3) .
\eeq
Therefore, proceeding as in Lemma ~\ref{Lem: Pprop}, we get~\eqref{eq: est 5.2}.
\end{proof}

\begin{proof}[Proof of Lemma~\ref{lem: gkpsi'}]
Choose $f\in  L^1_x(\mathbb{R}^n)$ and we have 
    \begin{equation}
    \begin{aligned}
        \|(1-F_{1,k})f\|_{L^2_x(\mathbb{R}^n)}&=\left(\int_{\mathbb{R}^n} |1-F_{1,k}(1+s+2^k)^{\frac14}|q|>1|^2|\hat{f}(q)|^2dq\right)^{\frac12}\\
        &\leq \left(\int_{\mathbb{R}^n} |1-F_{1,k}(1+s+2^k)^{\frac14}|q|>1|^2dq\right)^{\frac12}\|f\|_{L^1_x(\mathbb{R}^n)}\\
        &\lesssim_{n} \frac{1}{\langle1+s+2^k\rangle^{\frac n2}}\|f\|_{L^1_x(\mathbb{R}^n)}.
    \end{aligned}
    \end{equation}
This completes the proof.
\end{proof}

\noindent{\textbf{Acknowledgements}: }A. S. was partially supported by Simons Foundation Grant number 851844 and NSF grants DMS-2205931. X. W. was supported by Australian Laureate Fellowships, grant FL220100072. T. Z. was supported by National Natural Science Foundation of China (Grant No.11931010) and Natural Science Foundation of Zhejiang Province (Grant No. LDQ23A010001). This work was done by the second author while at the School of Mathematical Sciences, Zhejiang University and  South China University of Technology. Part of this work was completed by the third author while at Texas A\&M University and the Fields Institute.

	\vspace{2mm}
	
\noindent \textbf{Conflict of interest.} On behalf of all authors, the corresponding author states that there is no conflict of interest.
	\vspace{2mm}

	\noindent \textbf{Data availability statement.}
	Data sharing is not applicable to this article as no data sets were generated or analyzed during the current study.


\end{document}